\documentclass[11pt,leqno,openany]{article}
\usepackage{amsmath,amssymb,amsthm,textcomp,times,a4}
\usepackage[T1]{fontenc}
\usepackage[latin1]{inputenc}
\usepackage{colortbl}
\usepackage[english]{babel}
\usepackage{dsfont}

\newcommand{\bfm}[1]{\boldsymbol{#1}}

\newtheorem{theorem}[subsection]{Theorem}

\newtheorem{definition}[subsection]{Definition}
\newtheorem{lemma}[subsection]{Lemma}
\newtheorem{proposition}[subsection]{Proposition}

\newtheorem{assumption}[subsection]{Assumption}
\newtheorem{remark}[subsection]{Remark}

\def \D{ {\rm I}\!{\rm D} }
\def \F{ \mathbb  F}
\def \H{ \mathbb H}

\def \R{ \mathbb R }
\def \reel{ \mathbb R }
\def \nat{ \mathbb N}

\def \ent{ \mathbb Z}
\def \one{\mathbf{ 1}}%{\rm 1}\mkern-5mu{\rm I} }

\def \E{ {\rm I}\!{\rm E} }
\def \prob{ {\rm I}\!{\rm P} }
\def \P{ {\rm I}\!{\rm P} }

\def \tor{ \mathbb T}
\begin{document}
\author{Rémi Rhodes\footnote{Laboratoire d'Analyse Topologie Probabilités, Université de Provence, 39 rue Joliot Curie, 13453 Marseille Cedex 13, France,\hspace{2cm} e-mail: \textsf{rhodes@cmi.univ-mrs.fr}}}
\title{On homogenization of space-time dependent degenerate random flows}
\maketitle
\begin{abstract}
We study a diffusion with time-dependent random coefficients. The
diffusion coefficient is allowed to degenerate. We prove an
invariance principle when this diffusion is supposed to be
controlled by another one with time independent coefficients.
\end{abstract}

%%%%%%%%%%%%%%%%%%%%%%%%%%%%%%%%%%%%%%%%%%%%%%%%%%%%%%%%%%%%%%%%%%%%%%%%%%%%%%%%%%%%%%%%%%%%%%%%%%%%%%%%%

\section{Introduction}

%%%%%%%%%%%%%%%%%%%%%%%%%%%%%%%%%%%%%%%%%%%%%%%%%%%%%%%%%%%%%%%%%%%%%%%%%%%%%%%%%%%%%%%%%%%%%%%%%%%%%%%%%%%

We want to establish an invariance principle for a diffusive
particle in a random flow described by the following Stochastic
Differential Equation (SDE)
\begin{equation*}%\label{eqproceps}
X^{\omega }_t=x+\int_{0}^{t}b\left( r,X^{\omega}_r,\omega
\right)\,dr+\int_{0}^{t}\sigma \left(r,X^{ \omega}_r,\omega
\right)\,dB_r,
\end{equation*}
where $B$ is a d-dimensional Brownian motion and $\sigma,b$ are
stationary random fields. $b$ is defined in such a way that the
generator at time $t$ of the diffusion coincides on smooth functions
with
\begin{equation}\label{def_gen}
L^\omega=(1/2)e^{2V(x,\omega )}{\rm div}_x\big(e^{-2V(x,\omega
)}[a(t,x,\omega )+H(t,x,\omega )]\nabla _x\,\,\big)\,.
\end{equation}
Here $a(t,x,\omega  )$ is equal
 to $\sigma\sigma^*(t,x,\omega )$. $V$ and $H$ are stationary random fields, $V$ is bounded and $H$ antisymmetric.

\noindent We will then be in position to study the effective
diffusion on a macroscopic scale of the following
convection-diffusion equation
\begin{equation}\label{equationeps}
\begin{split}
\partial_t z(t,x,\omega ) =(1/2){\rm Trace}[a\triangle_{xx}z](t,x,\omega
)+b\cdot\nabla_xz(t,x,\omega),
\end{split}
\end{equation}
with certain initial condition. We will prove that, in probability
with respect to $\omega $,
 \begin{equation*} \lim_{\varepsilon
\rightarrow 0}z(t/\varepsilon ^2,x/\varepsilon ,\omega
)=\overline{z}(t,x)
\end{equation*}
where $\overline{z}$ is the solution of a deterministic equation
\begin{equation}\label{eq_lim}
\partial_t \overline{z}(t,x)={\rm
Trace}[A\triangle_{xx}\overline{z}](t,x) .
\end{equation} $A$ is a constant
matrix - the matrix of so-called effective coefficients.

Homogenization problems have been extensively studied in the case of
periodic flows (cf. \cite{bhat}, \cite{pardouxdeg},
\cite{pardouxper}, and many others). The study of random flows (see
\cite{olla}, \cite{osada}, \cite{sethu}, and many others) spread
rapidly thanks to the techniques of the {\it environment
 as seen from the particle} introduced by Kipnis and Varadhan in \cite{kipnis}, at least in the case of time independent
 random flows.  Recently, there have been results going beyond these techniques in the case of isotropic coefficients
 which are small perturbations of Brownian motion (see \cite{sznitman}).
But there are only a few works in the case of space-time dependent
random flows (see \cite{olla1} or \cite{olla2} for instance in the
case $\sigma={\rm Id}$). A quenched version of the invariance
principle is stated in \cite{fannjang} provided that the diffusion
coefficient satisfies a strong uniform non-degeneracy assumption. In
this case, the regularizing properties of the heat kernel are widely
used to face with the non-reversibility of the underlying processes.
Some results stated in Markovian flows are also established in
\cite{fannjang4} or \cite{fannjang5}. \\The novelty of this work
lies in the ergodic and regularizing properties required on the
coefficients, which are not far from being minimal. The only
restriction is the control of the diffusion process with an ergodic
and time independent one. As a consequence, this work includes the
static case where all the coefficients do not depend on time.
Moreover, these assumptions allow the diffusion matrix to
degenerate. Typically it can degenerate in certain directions or
vanish on subsets of null measure but cannot totally reduce to zero
on an open subset of $\reel^d$. However, considering such strong
degeneracies remains a quite open problem for random stationary
coefficients (for recent advances in the static periodic case, see
\cite{pardouxdeg}).

We will outline now the main ideas of the proof. Our goal is to show that the rescaled process
$$\varepsilon X^{\omega }_{t/\varepsilon ^2}=\varepsilon \int_{0}^{{t/\varepsilon ^2}}b\left( r,X^{\omega}_r,\omega
\right)\,ds+\varepsilon \int_{0}^{{t/\varepsilon ^2}}\sigma
\left(r,X^{ \omega}_r,\omega \right)\,dB_s$$ converges in law to a
Brownian motion with a certain positive covariance matrix. The
general strategy (see \cite{kozlov}) consists in finding an
approximation of the first term on the right-hand side by a family
of martingales and then in applying the central limit theorem for
martingales. To find such an approximation, we look at the
environment as seen from the particle
$$Y_t=\tau _{t,X^\omega _t}\omega ,$$ where $\{\tau _{t,x}\}$ is a group of measure preserving transformation on
a random medium $\Omega $ (see Definition \ref{medium}). Thanks to
the particular choice of the drift, an explicit invariant measure
can be found for this Markov process. The ergodicity is ensured by
the geometry of the diffusion coefficient $\sigma$ (see Assumptions
\ref{hypcontrol} and \ref{ergodicity}). The approximation that we
want to find leads to study the equation ($\lambda >0$)
\begin{equation}\label{eqresb}
   \lambda {\bfm u}_\lambda -({\bfm L}+D_t){\bfm u}_\lambda={\bfm b}
\end{equation}
 where ${\bfm L}+D_t$ coincides with the generator of the
process $Y$ on a certain class of functions (the term $D_t$ is due
to the time evolution and ${\bfm L}$ is an unbounded operator on the
medium $\Omega$ associated to \eqref{def_gen}). Here are arising the
difficulties resulting from the time dependence. Due to the term
$D_t$, the Dirichlet form associated to ${\bfm L}+D_t$ does not
satisfy any sector condition (even weak). However, for a suitable
function ${\bfm b}$, \eqref{eqresb} can be solved with the help of
an approximating sequence of Dirichlet forms with weak sector
condition. Then, usual techniques used in the static case fall short
of establishing the so-called sublinear growth of the
 correctors ${\bfm u}_\lambda$. To get round this difficulty, regularizing properties of the heat kernel are used in \cite{fannjang},
 \cite{olla1} or \cite{olla2}.
 Here the degeneracies of the diffusion coefficient prevents us from using such arguments. The strategy here consists in
  separating the time and spatial evolutions (see Assumption \ref{hypcontrol}). We introduce a new operator $\widetilde{{\bfm S}}$ whose
 coefficients do not depend on time. Then the spectral calculus linked to the normal operator $  \widetilde{{\bfm S}}+
 D_t$ will be determining to establish the desired estimates for the solution  ${\bfm v}_\lambda $ of the equation
$$  \lambda {\bfm v}_\lambda  -(\widetilde{{\bfm S}}+D_t){\bfm v}_\lambda   ={\bfm b}.$$
Finally, with perturbation methods, we show that these estimates
remain valid for the correctors $ {\bfm u}_\lambda $.

%%%%%%%%%%%%%%%%%%%%%%%%%%%%%%%%%%%%%%%%%%%%%%%%%%%%%%%%%%%%%%%%%%%%%%%%%%%%%%%%%%%%%%%%%

\section{Notations, Setup and Main Result}\label{notsetup}

%%%%%%%%%%%%%%%%%%%%%%%%%%%%%%%%%%%%%%%%%%%%%%%%%%%%%%%%%%%%%%%%%%%%%%%%%%%%%%%%%%%%%%%%%

Let us first introduce a random medium
\begin{definition}\label{medium}
Let $(\Omega ,{\cal G},\mu )$ be a probability space and
$\left\{\tau_{t,x};(t,x)\in \R\times\reel ^d\right\}$ a
stochastically continuous group of measure preserving
transformations acting ergodically on $\Omega $:

1) $\forall A\in {\cal G},\forall (t,x)\in \reel\times\reel^d$, $\mu
(\tau _{t,x}A)=\mu (A)$,

2) If for any $(t,x)\in \reel\times\reel^d$, $\tau _{t,x}A=A$ then
$\mu (A)=0$ or $1$,

3) For any measurable function ${\bfm g}$ on $(\Omega ,{\cal
G},\mu)$, the function $(t,x,\omega )\mapsto {\bfm g}(\tau_{t,x}
\omega)$ is measurable on $(\reel\times\reel^d\times\Omega ,{\cal
B}(\reel\times\reel^d)\otimes {\cal G})$.
\end{definition}
In what follows we will use the bold type to denote a function
${\bfm g}$ from $\Omega$ into $\R$ (or more generally into $\R^n$,
$n \geq 1$)  and the unbold type $g(t,x,\omega)$ to denote the
associated representation mapping $(t,x,\omega) \mapsto {\bfm
g}(\tau_{t,x} \omega)$. The space of square integrable functions on
$(\Omega ,{\cal G},\mu )$ is denoted by $L^2(\Omega)$, the usual
norm by $|\, \cdot \,|_2$ and the corresponding inner product by
$(\, \cdot\,,\, \cdot\,)_2$. Then, the operators on $L^2(\Omega )$
defined by $T_{t,x}{\bfm g}(\omega )={\bfm g}(\tau_{t,x}\omega )$
form a strongly continuous group of unitary maps in $L^2(\Omega )$.
Each function ${\bfm g}$ in $ L^2(\Omega )$ defines in this way a
stationary ergodic random field on $\reel ^{d+1}$. The group
possesses $d+1$ generators defined for $i=1,\dots,d,$ by
\begin{equation*} D_i{\bfm f}=\frac{ \partial }{\partial x_i}T_{0,x}{\bfm
f}|_{(t,x)=0},\quad  \text{ and }\quad D_t{\bfm f}=\frac{ \partial
}{\partial t}T_{t,0}{\bfm f}|_{(t,x)=0},
\end{equation*}
which are closed and densely defined. Denote by ${\cal C}$ the dense
subset of $L^2(\Omega )$ defined by
\begin{equation*}
{\cal C}=\mathrm{Span}\big\{{\bfm f}*\varphi ;{\bfm f}\in L^2(\Omega
),\varphi \in C^\infty _c(\reel ^{d+1} )\big\},\quad \text{ with }\
{\bfm f}*\varphi(\omega )=\int_{\reel ^{d+1}}{\bfm
f}(\tau_{t,x}\omega )\varphi (t,x)\,dt\,dx,
\end{equation*}
where $C^\infty _c(\reel ^{d+1} )$ is the set of smooth functions on
$\reel^{d+1}$ with a compact support. Remark that ${\cal C}\subset
{\rm Dom}(D_i)$ and $D_i({\bfm f}*\varphi)=-{\bfm f}* \frac{\partial
\varphi }{\partial x_i}$. This last quantity is also equal to
$D_i{\bfm f}*\varphi $ if ${\bfm f}\in {\rm Dom}(D_i)$.

Consider now the measurable functions ${\bfm
\sigma}:\Omega\rightarrow \reel^{d\times d}$, $\widetilde{{\bfm
\sigma}}:\Omega\rightarrow \reel^{d\times d}$, ${\bfm
H}:\Omega\rightarrow \reel^{d\times d}$ and ${\bfm
V}:\Omega\rightarrow \reel $ and assume that ${\bfm H}$ is
antisymmetric. Define ${\bfm a}={\bfm \sigma}{\bfm \sigma}^*$ and
$\widetilde{{\bfm a}}=\widetilde{{\bfm \sigma}}\widetilde{{\bfm
\sigma}}^*$. The function ${\bfm V}$ does not depend on time, that
means $\forall t\in\R$, $T_{t,0}{\bfm V}={\bfm V}$.

\begin{assumption}{{\bf (Regularity of the coefficients)}}\label{hypregularity}\\
$\bullet$ Assume that $\forall i,j,k,l=1,\dots,d$,  $\ {\bfm
a}_{ij},\widetilde{{\bfm a}}_{ij},{\bfm V}, {\bfm H}_{ij},D_l{\bfm
a}_{ij}\text{ and }D_l\widetilde{{\bfm a}}_{ij} \in {\rm
Dom}(D_k)$.\\
$\bullet$ Define, for $i=1,\dots,d$,
\begin{equation}
\begin{split}
\label{defb}
 {\bfm b}_i(\omega ) & =\sum_{j=1}^{d}\big(\frac{
1}{2}D_j{\bfm a}_{ij}(\omega)- {\bfm a}_{ij}D_j{\bfm V}(\omega
)+\frac{ 1}{2}
D_j{\bfm H}_{ij}(\omega )\big),\\
 \widetilde{{\bfm b}}_i(\omega ) &
=\sum_{j=1}^{d}\big(\frac{ 1}{2}D_j\widetilde{{\bfm
a}}_{ij}(\omega)- \widetilde{{\bfm a}}_{ij}D_j{\bfm V}(\omega
)\big),
\end{split}
\end{equation} and assume that the applications $(t,x)~\mapsto~b_i(t,x,\omega
)$, $(t,x)~\mapsto~\widetilde{b}_i(t,x,\omega )$,
$(t,x)\mapsto\sigma(t,x,\omega )$
%, $x~\mapsto~\widetilde{b}(t,x,\omega)$ and $x\mapsto\widetilde{\sigma}(t,x,\omega )$
are globally Lipschitz. Moreover, the coefficients ${\bfm \sigma}$,
${\bfm a}$, ${\bfm b}$, $\widetilde{{\bfm \sigma}}$, ${\bfm V}$,
${\bfm H}$ are uniformly bounded by a constant $K$. (In particular,
this ensures existence and uniqueness of a global solution of SDE
\eqref{diffusion}.)
\end{assumption}

Here is the main assumption of this paper
\begin{assumption}{{\bf (Control of the
coefficients)}}\label{hypcontrol}\\
$\bullet$ $\widetilde{{\bfm \sigma  }}$ does not depend on time
(i.e. $\forall t\in \reel$, $T_t\widetilde{{\bfm \sigma
}}=\widetilde{{\bfm \sigma  } }$) and ${\bfm H},{\bfm a}\in {\rm Dom}(D_t)$.
As a consequence, the matrix $\widetilde{{\bfm a}}$ does not depend on time either.\\
$\bullet$  There exist five positive constants $m, M, C_1^H, C^H_2,
 C^a_2$ such that, $\mu$ a.s.,
\begin{equation}\label{inega}
m\widetilde{{\bfm a}}\leq {\bfm a}\leq M\widetilde{{\bfm a}},
\end{equation}
\begin{equation}\label{matrixh}
|{\bfm H}|\leq C^H_1\widetilde{{\bfm a}},\quad |D_t{\bfm H}|\leq
C^H_2\widetilde{{\bfm a}} \quad \text{and} \quad |D_t{\bfm a}|\leq
C^a_2\widetilde{{\bfm a}},
\end{equation}
where $|{\bfm A}|$ stands for the symmetric positive square root of
${\bfm A}$, i.e. $|{\bfm A}|=\sqrt[]{ -{\bfm A}^2}$.
\end{assumption}
\noindent For instance, if the matrix ${\bfm a}$ is uniformly
elliptic and bounded, $\widetilde{{\bfm \sigma} }$ can be chosen as
equal to the identity matrix ${\bfm {\rm Id}} $ and then
$\eqref{matrixh} \Leftrightarrow {\bfm H},\,D_t{\bfm H}\text{ and
}D_t{\bfm a}\in L^\infty (\Omega )$.\\ Let us now set out the
ergodic properties of this framework
\begin{assumption}\label{ergodicity}{{\bf (Ergodicity) }}
Let us consider the operator $\widetilde{{\bfm S}}= (1/2)e^{ 2{\bfm
V}}\sum_{i,j=1}^{d}D_i(e^{- 2{\bfm V}}\widetilde{{\bfm a}}_{ij}D_j \
)$ with domain ${\cal C}$. From Assumption \ref{hypregularity}, we
can consider its Friedrich extension (see \cite[Ch. 3, Sect.
3]{fukushima}) which is still denoted $\widetilde{{\bfm S}}$. Assume
that each function ${\bfm f}\in {\rm Dom}(\widetilde{{\bfm S}})$
satisfying $\widetilde{{\bfm S}}{\bfm f}=0$ must be $\mu$ almost
surely equal to some function that is invariant under space
translations.
\end{assumption}

Even if it means adding to ${\bfm V}$ a constant (and this does not
change the drift ${\bfm b}$, see \eqref{defb}), we make the
assumption that $\int e^{-2{\bfm V}}\,d\mu =1$. Thus we can define a
new probability measure on $\Omega $ by
\begin{equation*}
d\pi (\omega )=e^{-2{\bfm V}(\omega )}\,d\mu (\omega ).
\end{equation*}

We now consider a standard $d$-dimensional Brownian motion defined
on a probability space $(\Omega ',{\cal F},\prob)$ (the medium and
the Brownian motion are mutually independent) and the diffusions in
random medium given as the solutions of the following Stochastic
Differential Equations (SDE)
\begin{equation}\label{diffusion}
\begin{split}
X^{\omega }_t & =x+\int_{0}^{t}b\left(r, X^{\omega }_r,\omega
\right)\,dr+ \int_{0}^{t}\sigma\left(r,X^{\omega }_r,\omega
\right)\,dB_r,\\
 \widetilde{X}^{\omega}_t & =x+\int_0^t\widetilde{b}(X^{\omega }_r,\omega )\,dr
+\int_0^t \widetilde{\sigma}(X^{\omega }_r,\omega ) \,dB_r.
\end{split}
\end{equation}
The main result of this paper is stated as follows
\begin{theorem}
The law of the rescaled process $\varepsilon
X^\omega_{t/\varepsilon^2}$  converges in probability (with respect
to $\omega$) to the law of a Brownian motion with a certain
covariance matrix A (see \eqref{def_A}).
\end{theorem}

%%%%%%%%%%%%%%%%%%%%%%%%%%%%%%%%%%%%%%%%%%%%%%%%%%%%%%%%%%%%%%%%%%%%%%%%%%%%%%%%%%%%%%%%%%

\section{Examples}

%%%%%%%%%%%%%%%%%%%%%%%%%%%%%%%%%%%%%%%%%%%%%%%%%%%%%%%%%%%%%%%%%%%%%%%%%%%%%%%%%%%%%%%%%%

There are many ways to ensure the validity of Assumption
\eqref{ergodicity}. In particular, it is satisfied when, for almost
all $\omega \in \Omega $, the $\reel ^d$-valued Markov process
$\widetilde{X}^{\omega} $, whose generator coincides on smooth
functions with
$$\widetilde{{\bfm S}}^{\omega}=\frac{e^{2V(x,\omega )}}{2}{\rm Div}_x\left(e^{-2V(x,\omega )}\widetilde{a}(x,\omega )
\nabla_x\,.\right), $$  is irreducible in the following sense.
Suppose that, starting from any point of $\reel ^d$, the process
reaches each subset of $\reel ^d$ of non-null Lebesgue measure in
finite time. That means that there exists a measurable subset
$N\subset \Omega $ with $\mu (N)=0$ such that $\forall \omega \in
\Omega \setminus N$, for each measurable subset $B$ of $\reel^d$
with $\lambda _{Leb}(B)>0$, $\forall x\in\reel^d,\exists t>0,$
\begin{equation}\label{eq_probpositive}
\prob_x\left(\widetilde{X}^{\omega }_t\in B\right)>0.
\end{equation}
This can be proved as in \cite{olla2} section 3 or in \cite{olla}
chapter 2 Theorem 2.1, in studying the $\Omega$-valued Markov
process
$\widetilde{Y}_t(\omega)=\tau_{0,\widetilde{X}^\omega_t}\omega$,
whose generator coincides on ${\cal C}$ with $\widetilde{{\bfm S}}$.
As an easy consequence, if the diffusion coefficient
$\widetilde{{\bfm a}}$ is uniformly elliptic or satisfies a strong
Hörmander condition (see \cite{kusuoka} for further details), then
estimates on the transition densities of the process
$\widetilde{X}^\omega$ ensure \eqref{eq_probpositive}.

Let us now tackle the issue of constructing examples that do not
satisfy any uniform ellipticity assumption or even strong Hörmander
condition. In what follows, two examples are given. The first one
deals with periodic coefficients. The second one is a random medium
with a random chessboard structure and thereby does not reduce to
the periodic case.

\subsection{A periodic example}
%%%%%%%%%%%%%%%%%%%%%%%%%%%%%%%%%%%%%%%%%%%%%%%%%%%%%%%%%%%%%%%%%%%%%%%%%%%%%%%%%%%%%%%%%%%%%%%%%%%%%%%%%%%%%%%%%%%%%%%%

Let us construct a periodic example on the torus $\tor^3$, where the
diffusion matrix reduces to zero on a certain subset with null
Lebesgue measure. We define a time-independent matrix-valued
function
$$\widetilde{\bfm{\sigma}}(t,x,y)=(1-\cos(x))(1-\cos(y))\left(\begin{array}{cc}
 1 &  0 \\
 0 &  1 \end{array}\right).$$
For simplicity, we choose ${\bfm V}=\bfm{H}(x,y)=0$. Thanks to the
(not uniform!) ellipticity of the diffusion coefficient inside the
cell $\verb"C"=]0,2\pi[\times]0,2\pi[$, it is not very difficult to
see that \eqref{eq_probpositive} is satisfied. Indeed, each subset
$B\subset [0;2 \pi]^2$ with a strictly positive Lebesgue measure
necessarily satisfies $ \lambda _{Leb}(B\cap \verb"C")>0$. As
explained above, this is sufficient to ensure Assumption
\ref{ergodicity}. Let us now focus on Assumption \ref{hypcontrol}.
The strategy consists in choosing a smooth function ${\bfm U}:
\tor^3\rightarrow \R^{2\times 2}$ satisfying $\alpha^{-1} {\rm
Id}\leq {\bfm U}{\bfm U}^*(t,x,y)\leq \alpha {\rm Id}$ for some
$\alpha>0$, and then in defining ${\bfm
\sigma}(t,x,y)=\widetilde{\bfm{\sigma}}(t,x,y)\bfm{U}(t,x,y)$, for
which Assumption \ref{hypcontrol} is easily checked.

\subsection{An example on chessboard structures}
%%%%%%%%%%%%%%%%%%%%%%%%%%%%%%%%%%%%%%%%%%%%%%%%%%%%%%%%%%%%%%%%%%%%%%%%%%%%%%%%%%%%%%%%%%%%%%%%%%%%%%%%%%%%%%%%%%%%%%%%%%%%%

Let us now explain how to construct a random medium with chessboard
structures. Given $d\geq 1$, consider a sequence $(\varepsilon
_{(k_1,\dots,k_d)})_{(k_1,\dots,k_d)\in \ent ^d}$ of independant
Bernouilli random variables with parameter $p\in ]0,1[$ and define a
process $\widetilde{\eta }$ as follows: for each $x\in \reel ^d$,
there exists a unique $(k_1,\dots,k_d)\in \ent  ^d$ such that $x$
belongs to the cube $[k_1,k_1+1[\times \dots\times [k_d,k_d+1[$.
Then define the process $\widetilde{\eta}:\reel^d\rightarrow\R$ by:
$\forall x\in\reel^d, \quad \widetilde{\eta}_x= \varepsilon
_{(k_1,\dots,k_d)}$. The law of this process is invariant and
ergodic with respect to $\ent ^d$ translations. Roughly speaking, we
are drawing a $d$-dimensional chessboard on $\reel^d$, for which we
are coloring each cube of the chessboard either in black with
probability $p$ or in white with probability $1-p$. It remains to
make the process invariant under $\reel ^d$ translations. To this
purpose, choose a uniform variable $U$ on the cube $[0,1[^d$
independent of the sequence
$(\varepsilon_{(k_1,\dots,k_d)})_{(k_1,\dots,k_d)\in \ent ^d}$ and
define for $x\in\R^d$, $\overline{\eta }_x=\widetilde{\eta }_{x+U}$.
In a way, this corresponds to a random change of the origin of the
chessboard. It can be checked that we get a stationary ergodic
random field on $\reel ^d$. Let us now tackle the issue of the
regularity of the trajectories. Consider a $C^\infty (\reel ^d)$
 function $\varphi $ with a compact and very small support (for instance, included in the ball $B(0,1/4)$) and define a new process
$\eta _x=\int_{\reel ^d}\overline{\eta }_y\varphi
(x-y)\,dy=\overline{\eta }*\varphi (x)$, which is a stationary
ergodic random process with smooth trajectories. That is enough for
a general framework.

Let us now consider the process $\omega
_{(t,x)}=(\beta_t,\alpha^1_{x_1},\alpha^2_{x_2})_{t\in \reel
,x=(x_1,x_2)\in \reel ^2}$, where the three processes
$\alpha^1,\alpha ^2 $ and $\beta $ are mutually independent and
constructed as prescribed above. Hence \sloppy $\left\{\omega
_{(t,x)};(t,x)\in\reel\times\reel^2\right\}$ is an ergodic
stationary process and we can consider the random medium
$\Omega=C(\reel\times\reel ^{2};\reel ^3)$ equipped with the
probability law of this process. \\We define the matrix
$\widetilde{{\bfm \sigma}}(\omega  )=\left[\begin{array}{cc}
1 & 0\\
0 & \alpha^1_{0}
\end{array}\right]$
and ${\bfm V}=0$ (or any bounded function of the random field
$\alpha $). We can choose any matrix-valued function ${\bfm
U}:\Omega \rightarrow \R^{2\times 2}$ such that ${\bfm U}{\bfm U}^*$
is uniformly elliptic and bounded, and then we set ${\bfm \sigma
}=\widetilde{{\bfm \sigma }}{\bfm U}$. It can be proved that
Assumption \ref{ergodicity} is satisfied. Actually, the ergodicity
property for $\tilde{{\bfm \sigma}}$ is very intuitive. Indeed, the
matrix $\tilde{\sigma}(\cdot,\omega)$ degenerates only on some
stripes (the white ones), and in fact only on a part of each of them
(depending on the support of $\varphi$), and only along the
$y_2$-axis direction: while lying on the degenerating part of a
white stripe, the diffusion associated to $(1/2) \sum_{i,j=1}^2
\partial_i( \tilde{a} _{i,j} \partial_j)$ can
only move along the $y_1$-axis direction. Nevertheless, with
probability $1$, the process encounters a black stripe sooner or
later (because the parameter $p$ belongs to $]0,1[$): it thus
manages to move up and down and to reach every subset of the space.
 Ergodicity follows. Rigorous
arguments are however left to the reader.

%We just outline the main ideas. Let $(\varepsilon _k)_{k\in\ent}$ be
%the independant Bernouilli random variables of parameter $p\in]0;1[$
%used to
% construct the process $\alpha^1 $. There exists $N\subset \Omega $ such that $\mu (N)=0$ and $\forall \omega \in
% \Omega \setminus N$, there exists $k\in \ent$ such that $\varepsilon _k=1$. Hence, by the construction of $\alpha^1 $
% (and more precisely the fact that ${\rm Supp}(\varphi )~\subset~B(0,1/4)$ ), there exists $z\in \reel$ such that
% $\alpha^1(x,\omega ) =1$ for $x\in B_z=\{(x,y)\in \reel^2;|x-z|<  \frac{1}{4}\}$.
%
%Then, starting from any point of the plan, it is not difficult to
%see that the process $X^{\omega ,\widetilde{S}}$ can reach the strip
%$B_z$, because of the non-degeneracy of the diffusion coefficient
%along the x-direction. Then the reversibility of this process
%ensures that, starting from $B_z$, we can reach each subset of
%$\reel^2$. Finally, by the uniform ellipticity of the diffusion
%coefficient over $B_z$ and the strong Markov property, we can prove
%that the assumption of Proposition \ref{irred} is valid. The
%ergodicity is proved.

We can also consider a non-null stream matrix ${\bfm H}$. For
instance the matrix-valued function ${\bfm H}(\omega
)=\left[\begin{array}{cc}
0 & (\alpha^1_0) ^2\beta_0 \\
-(\alpha_0^1) ^2\beta_0  & 0
\end{array}\right],$ fits Assumption \ref{hypcontrol}.

%%%%%%%%%%%%%%%%%%%%%%%%%%%%%%%%%%%%%%%%%%%%%%%%%%%%%%%%%%%%%%%%%%%%%%%%%%%%%%%%%%%%%%%%%%

\section{Environment as seen from the particle}\label{sec_env}

%%%%%%%%%%%%%%%%%%%%%%%%%%%%%%%%%%%%%%%%%%%%%%%%%%%%%%%%%%%%%%%%%%%%%%%%%%%%%%%%%%%%%%%%%%

We now look at the {\it environments as seen from the particle}
associated to the processes $X$ and $\widetilde{X}$: they both are
$\Omega$-valued Markov processes and are defined by
\begin{equation}
\widetilde{Y}_t(\omega )  =\tau_{t,\widetilde{X}^{\omega }_t}\omega
,\quad\text{ and }\quad Y_t(\omega )  =\tau_{t,X^{\omega }_t}\omega
,
\end{equation}
where the processes  $X^{\omega }$ and $\widetilde{X}^{\omega }$
both starts from the point $0\in \reel ^d$. An easy computation
proves that the generators of these Markov processes respectively
coincide on ${\cal C}$ with $\widetilde{{\bfm S}}+D_t$ and ${\bfm
L}+D_t$, where ${\bfm L}$ is defined on ${\cal C}$ by
\begin{equation}\label{defL}
%\begin{split}
%\widetilde{{\bfm S}}+D_t & =\frac{e^{2{\bfm V}}}{2}\sum_{i,j=1}^dD_i\big(e^{-2{\bfm V}}\widetilde{{\bfm a}}_{ij}D_j \,)+D_t,\\
{\bfm L} = \frac{e^{2{\bfm V}}}{2}\sum_{i,j=1}^dD_i\big(e^{-2{\bfm
V}}[{\bfm a}+{\bfm H}]_{ij}D_j \,).
%\end{split}
\end{equation}
Hence $\pi $ is an invariant measure for both processes (see also
\cite{oelschlager}). Both associated semigroups thus extend
continuously to $L^2(\Omega,\pi)$. We should point out that the
invariant measure need not be unique.

%%%%%%%%%%%%%%%%%%%%%%%%%%%%%%%%%%%%%%%%%%%%%%%%%%%%%%%%%%%%%%%%%%%%%%%%%%%%%%%%%%%%%%%%%%

\section{Poisson's equation}\label{poisson}

%%%%%%%%%%%%%%%%%%%%%%%%%%%%%%%%%%%%%%%%%%%%%%%%%%%%%%%%%%%%%%%%%%%%%%%%%%%%%%%%%%%%%%%%%%

The aim of this section is, at first, to find a solution ${\bfm
u}_\lambda$ of the resolvent equation that can formally be rewritten
(a rigorous definition of each term is given later), for $\lambda
>0$, as:
\begin{equation}\label{eqpoisson}
\lambda {\bfm u}_\lambda -({\bfm L}+D_t){\bfm u}_\lambda ={\bfm h}.
\end{equation}
Since the associated Dirichlet form satisfies no sector condition
(even weak), existence and regularity of such a solution is
generally a tricky work, especially in considering degeneracies both
in time and in space. However, for a suitable right-hand side, this
equation can be solved with the help of an approximating sequence of
Dirichlet forms satisfying a weak sector condition. Thereafter we
study the asymptotic behaviour of the solution ${\bfm u}_\lambda$ as
$\lambda\rightarrow 0$.

\subsection{Setup}\label{setup}
%%%%%%%%%%%%%%%%%%%%%%%%%%%%%%%%%%%%%%%%%%%%%%%%%%%%%
Let us denote by $(\widetilde{P}_t)_t$ the semigroup on
$L^2(\Omega,\pi )$ generated by the process $\widetilde{Y}$ and by
$(\widetilde{P}^*_t)_t$ its adjoint operator. Let us also denote by
$(\bar{P}_t)_t$ the self-adjoint semigroup on $L^2(\Omega,\pi )$
generated by the process
$\bar{Y}_t(\omega)=\tau_{0,\widetilde{X}^{\omega }_t}\omega$. Its
generator is $\widetilde{{\bfm S}} $. From the time independence of
the coefficients $\widetilde{{\bfm b}}$ and $\widetilde{{\bfm
\sigma}}$, it is readily seen that, that $\forall {\bfm f}\in
L^2(\Omega ,\pi )$, $\widetilde{P}_t{\bfm f}=T_{t,0}\bar{P}_t{\bfm
f}=\bar{P}_tT_{t,0}{\bfm f}$. As a consequence,
$\widetilde{P}^*_t=T_{-t,0}\bar{P}_t{\bfm f}=\bar{P}_tT_{-t,0}{\bfm
f}$, in such a way that
$$\widetilde{P}_t(\widetilde{P}_t^*{\bfm f})=\widetilde{P}_t^*(\widetilde{P}_t{\bfm f}).$$
The generator in $L^2(\Omega,\pi)$ of $(\widetilde{P}_t)_t$, wrongly
denoted by $[\widetilde{{\bfm S}} +D_t]$, is then normal (see
Theorem 13.38 in \cite{revuz}) so that we can find a spectral
resolution of the identity $E$ on the Borelian subsets of
$\reel_+\times \reel $ such that
$$-[\widetilde{{\bfm S}} +D_t]=\int_{\reel_+\times \reel}(x+iy)\,E(dx,dy).$$
 Actually, we have $-\widetilde{{\bfm S}}=\int_{\reel_+\times \reel}x\,E(dx,dy),\quad \text{ and
 }\quad-D_t=\int_{\reel_+\times\reel}iy\,E(dx,dy)$.
Indeed, $\widetilde{{\bfm S}}$ and $\int_{\reel_+\times
\reel}x\,E(dx,dy) $ are both self-adjoint and coincide on ${\cal
C}$. From \cite[Ch. 1, Sect. 3]{fukushima}, they are equal. The same
arguments hold for $D_t$ and
$\int_{\reel_+\times\reel}iy\,E(dx,dy)$. \\For any ${\bfm \varphi
},{\bfm \psi }\in L^2(\Omega )$, denote by $E_{{\bfm \varphi },{\bfm
\psi }}$ the measure defined by $E_{{\bfm \varphi },{\bfm \psi
}}=(E{\bfm \varphi },{\bfm \psi })_2$. From now on, denote by
$(.\,,.\,)_2$ the usual inner product in $L^2(\Omega,\pi )$. For any
${\bfm \varphi },{\bfm \psi }\in {\cal C}$, define
\begin{equation}
\langle{\bfm \varphi },{\bfm \psi }\rangle_{1}=\int_{\reel_+\times
\reel}x\,E_{{\bfm \varphi },{\bfm \psi }}(dx,dy)=-( {\bfm \varphi
},\widetilde{{\bfm S}}{\bfm \psi })_2
\end{equation}
and $\|{\bfm \varphi }\|_{1}=\sqrt[]{\langle{\bfm \varphi },{\bfm
\varphi }\rangle_{1} }$. By virtue of Assumption \eqref{inega}, this
semi-norm is equivalent on ${\cal C}$ to the semi-norm defined by
$\sqrt{-({\bfm \varphi },{\bfm S}{\bfm \varphi })_2}$,
\begin{equation}\label{equivnorm}
 m\|{\bfm \varphi }\|^2_{1}\leq -({\bfm \varphi },{\bfm S}{\bfm \varphi })_2\leq M\|{\bfm \varphi
}\|^2_{1},
\end{equation}
 where ${\bfm S}$ is the Friedrich extension of the
operator defined on ${\cal C}$ by $(1/2)e^{2{\bfm
V}}\sum_{i,j}D_i\big(e^{-2{\bfm V}}{\bfm a}_{ij}D_j\,\big) $.

Let $\F$ (respectively $\H$) be the Hilbert space that is the
closure of ${\cal C}$ in $L^2(\Omega )$ with respect to the inner
product $\varepsilon $ (resp. $\kappa $) defined on ${\cal C}$ by
\begin{equation*}
\begin{split}
&\varepsilon ({\bfm \varphi },{\bfm \psi })  =({\bfm \varphi },{\bfm
\psi })_2+ \langle{\bfm \varphi },{\bfm \psi }\rangle_{1}+(D_t{\bfm
\varphi
},D_t{\bfm \psi })_2\\
 & (\text{resp. }\kappa ({\bfm \varphi },{\bfm \psi })  =({\bfm
\varphi },{\bfm \psi })_2+ \langle{\bfm \varphi },{\bfm \psi
}\rangle_{1}).
\end{split}
\end{equation*}
Define the space $\D$ as the closure in $(L^2(\Omega ),|\,.\,|_2)$
of the subspace $\{(-\widetilde{S})^{1/2}{\bfm \varphi };{\bfm
\varphi }\in {\cal C}\}$. For any ${\bfm \varphi }\in {\cal C} $,
define $\Phi\big((-\widetilde{S})^{1/2}{\bfm \varphi }\big)={\bfm
\sigma}^*D_x{\bfm \varphi }\in (L^2(\Omega))^d$ and note that
$|\Phi\big((-\widetilde{S})^{1/2}{\bfm \varphi }\big)|_2^2=-({\bfm
\varphi },{\bfm S}{\bfm \varphi })_2 $. From \eqref{equivnorm},
$\Phi$ can be extended to the whole space $\D$ and this extension is
a linear isomorphism from $\D$ into a closed subset of
$(L^2(\Omega))^d$. Hence, for each function ${\bfm u}\in \H$, we define $\nabla^\sigma{\bfm u}=\Phi((-\widetilde{S})^{1/2}{\bfm u})$
 and this stands, in a way, for the gradient of ${\bfm u}$ along the direction $\sigma$.\\
 For each ${\bfm f}\in L^2(\Omega )$ satisfying
$\int_{\reel_+\times \reel}\frac{1}{x}\,E_{{\bfm f },{\bfm f
}}(dx,dy)<\infty $, we define
\begin{equation}
\|{\bfm f}\|^2_{-1}=\int_{\reel_+\times \reel}\frac{1}{x}\,E_{{\bfm f },{\bfm f }}(dx,dy).
\end{equation}
We point out that $\|{\bfm f}\|_{-1}<\infty $ if and only if there exists $C\in\reel$ such that for any  ${\bfm \varphi }\in {\cal C}$,
 $({\bfm f},{\bfm \varphi })_2\leq C\|{\bfm \varphi
 }\|_1$. For such a function ${\bfm f}$, $\|{\bfm f}\|_{-1}$ also
 matches the smallest $C$ satisfying this inequality. Remark that $\|{\bfm f}\|_{-1}<\infty
 $ implies $\pi({\bfm f})=0$.
%$$\inf\left\{C\in \reel;\forall {\bfm \varphi }\in {\cal C},\int_\Omega {\bfm f}{\bfm \varphi }\,d\pi \leq C\|{\bfm \varphi }\|_1\right\}=\|{\bfm f}\|_{-1}.$$
Denote by $ \H_{-1}$  the closure of $L^2(\Omega)$ in $\H^*$
(topological dual of $\H$) with respect to the norm $\|\
\|_{-1}$.\\
Let us now focus on the antisymmetric part ${\bfm H}$. We have
\begin{equation}\label{inegh}
|(u,{\bfm H}v)| \leq (u,|{\bfm H}|u)^{1/2}(v,|{\bfm H}|v)^{1/2}\leq
C^H_1(u,\widetilde{\bfm a}u)^{1/2} (v,\widetilde{{\bfm a}}v)^{1/2}.
\end{equation}
The second inequality follows from (\ref{matrixh}) and the first one
is a general fact of linear algebra.
%which can be shown by reducing the antisymmetric matrix ${\bfm H}$ in orthogonal basis to blocks of the type $\left(\begin{array}{cc} 0 & \alpha \\
%-\alpha & 0\end{array}\right)$. Then we just have to verify this inequality for such a type of matrix and it is a quite simple exercise.
We deduce
$$\forall {\bfm \varphi },{\bfm \psi }\in {\cal C},\quad (1/2)({\bfm H}D_x{\bfm \varphi }, D_x{\bfm \psi })_2 \leq
C_1^H \|{\bfm \psi }\|_{1}\|{\bfm \varphi  }\|_{1}.$$ Thus there
exists an antisymmetric  continuous bilinear form ${\bfm T}_H$ on
$\D\times \D $ such that
\begin{equation}\label{def_T}
\forall {\bfm \varphi },{\bfm \psi  }\in {\cal C},\quad (1/2)({\bfm
H}D_x{\bfm \varphi }, D_x{\bfm \psi })_2 ={\bfm
T}_H\big((-\widetilde{S})^{1/2}
 {\bfm \varphi },(-\widetilde{S})^{1/2}{\bfm \psi  }\big).
\end{equation}
Likewise, with the help of Assumption \ref{hypcontrol}, we define
the continuous bilinear forms ${\bfm T}_a$, $\partial_t{\bfm T}_a$,
$\partial_t{\bfm T}_H$, $\Lambda_s{\bfm T}_a$, $\Lambda_s{\bfm T}_a$
on $ \D\times \D \subset L^2(\Omega ,\pi )\times L^2(\Omega ,\pi )$
as follows: $\forall {\bfm \varphi },{\bfm \psi  }\in {\cal C}$,
\begin{align}
 (1/2)({\bfm a}D_x{\bfm \varphi }, D_x{\bfm \psi })_2 & ={\bfm
T}_a\big((-\widetilde{S})^{1/2}
 {\bfm \varphi },(-\widetilde{S})^{1/2}{\bfm \psi  }\big),\nonumber\\
  (1/2)(D_t{\bfm a}D_x{\bfm \varphi }, D_x{\bfm \psi })_2 &
=\partial_t{\bfm T}_a\big((-\widetilde{S})^{1/2}
 {\bfm \varphi },(-\widetilde{S})^{1/2}{\bfm \psi  }\big),\nonumber\\
 (1/2)(D_t{\bfm H}D_x{\bfm \varphi }, D_x{\bfm \psi })_2 &
=\partial_t{\bfm T}_H\big((-\widetilde{S})^{1/2}
 {\bfm \varphi },(-\widetilde{S})^{1/2}{\bfm \psi  }\big),\nonumber\\
(1/2)(\Lambda_s{\bfm a}D_x{\bfm \varphi }, D_x{\bfm \psi })_2 &
=\Lambda_s{\bfm T}_a\big((-\widetilde{S})^{1/2}
 {\bfm \varphi },(-\widetilde{S})^{1/2}{\bfm \psi  }\big),\nonumber\\
 (1/2)(\Lambda_s{\bfm H}D_x{\bfm \varphi }, D_x{\bfm \psi })_2 &
=\Lambda_s{\bfm T}_H\big((-\widetilde{S})^{1/2}
 {\bfm \varphi },(-\widetilde{S})^{1/2}{\bfm \psi  }\big),\nonumber
\end{align}
where, for any $s\in\R^*$, $\Lambda_s$ denotes the $L^2$-continuous
difference operator (remind of the definition of $T_{s,0}$ in
section \ref{notsetup}):
\begin{equation}\label{def_diffop}
 \forall {\bfm f}\in
L^2(\Omega),\quad \Lambda_s({\bfm f})=(T_{s,0}{\bfm f}-{\bfm f})/s.
\end{equation}
From Assumption \ref{hypcontrol}, the norms of the forms
$\Lambda_s{\bfm T}_a$ and $\Lambda_s{\bfm T}_H$ are uniformly
bounded with respect to $s\in\R^*$ and the forms are weakly
convergent respectively towards $\partial_t{\bfm T}_a$ and
$\partial_t{\bfm T}_H$.\\
Now, denote by ${\cal H}$ the subspace of $\H_{-1}$ whose elements
satisfy the condition: $\exists C>0,\forall s>0$ and $\forall {\bfm
\varphi }\in {\cal C}$, $ \langle{\bfm h} ,\Lambda_s{\bfm \varphi
}\rangle_{-1,1}\leq C \| {\bfm \varphi }\|_1$. For any ${\bfm h}\in
{\cal H}$, the smallest $C$ that satisfies such a condition is
denoted $\|{\bfm h}\|_T$. Then ${\cal H}$ is closed for the norm
$\|\,\|_{\cal H}=\|\,\|_{-1}+\|\,\|_T$. \\
Finally, let us now extend the operator ${\bfm L}$ defined on ${\cal
C}$ by \eqref{defL}. For any $\lambda>0$, consider the continuous
bilinear form ${\cal B}_\lambda$ on $\H\times\H$ that coincides on
${\cal C}\times{\cal C}$ with
\begin{equation*}%\label{}
\forall{\bfm \varphi},{\bfm \psi}\in {\cal C},\quad {\cal
B}_\lambda({\bfm \varphi},{\bfm \psi})=\lambda({\bfm \varphi},{\bfm
\psi})_2+[{\bfm T}_a+{\bfm T}_H]\big((-\widetilde{S})^{1/2}
 {\bfm \varphi },(-\widetilde{S})^{1/2}{\bfm \psi  }\big).
\end{equation*}
Thanks to Assumption \ref{hypcontrol} and the antisymmetry of ${\bfm
H}$, this form is clearly coercive. Thus it defines a strongly
continuous resolvent operator and consequently, the generator ${\bfm
L}$ associated to this resolvent operator. More precisely, ${\bfm
\varphi}\in\H$ belongs to ${\rm Dom}({\bfm L})$ if and only if
${\cal B}_\lambda({\bfm \varphi},\cdot)$ is $L^2$-continuous. In
this case, there exists ${\bfm f}\in L^2(\Omega)$ such that ${\cal
B}_\lambda({\bfm \varphi},\cdot)=({\bfm f},\cdot)_2$ and ${\bfm
L}{\bfm \varphi}$ is equal to ${\bfm f}-\lambda{\bfm \varphi}$. It
can be proved that this definition is independent of $\lambda>0$
(see \cite[Ch. 1, Sect. 2]{ma} for further details). Let us
additionally mention that the adjoint operator ${\bfm L}^*$ of
${\bfm L}$ in $L^2(\Omega,\pi)$ can also be described through ${\cal
B}_\lambda$. Indeed, ${\rm Dom}({\bfm L}^*)=\{{\bfm
\varphi}\in\H;{\cal B}_\lambda(\cdot,{\bfm \varphi})\text{ is
}L^2(\Omega)\text{-continuous.}\}$. If ${\bfm \varphi}\in{\rm
Dom}({\bfm L}^*)$, there exists ${\bfm f}\in L^2(\Omega)$ such that
${\cal B}_\lambda(\cdot,{\bfm \varphi})=({\bfm f},\cdot)_2$ and
${\bfm L}^*{\bfm \varphi}$ is equal to ${\bfm f}-\lambda{\bfm
\varphi}$.
\begin{remark}
For each function ${\bfm \varphi}\in {\cal C}\subset \H$, the
application ${\bfm L}{\bfm \varphi}$ can be viewed as a function of
$\H_{-1}$. Indeed, $\forall {\bfm \psi}\in {\cal C}$, $({\bfm
L}{\bfm \varphi},{\bfm \psi})_2=-[{\bfm T}_a+{\bfm
T}_H]\big((-\widetilde{{\bfm S}})^{1/2}{\bfm
\varphi},(-\widetilde{{\bfm S}})^{1/2}{\bfm \psi}\big)\leq
[M+C^H_1]\|{\bfm \varphi}\|_1\|{\bfm \psi}\|_1$. Hence, the
application ${\bfm \varphi}\mapsto{\bfm L}{\bfm \varphi}\in \H_{-1}$
can be extended to the whole space $\H$ so that, for each function
${\bfm u}\in\H$, we can define ${\bfm L}{\bfm u}$ as an element of
$\H_{-1}$ even if ${\bfm u}\not\in{\rm Dom}({\bfm L})$.
\end{remark}
\subsection{Existence of a solution:}\label{existence}
%%%%%%%%%%%%%%%%%%%%%%%%%%%%%%%%%%%%%%%%%%%%%%%%%%%%%%%%%%%%%%%%%%
This section is devoted to proving existence of solutions of
equation \eqref{eqpoisson} for a suitable right-hand side. The
difficulty lies in the strong degeneracy of the associated Dirichlet
form.
It satisfies no sector condition, even weak. However, it can be approximated by a family of Dirichlet forms with weak sector condition.\\
For any $\theta\in\{0;1\}$, $\lambda>0$ and $\delta\geq 0$, define
$B^\theta_{\lambda,\delta}$ as the (non-symmetric) bilinear
continuous form on $\F\times \F$ that coincides on ${\cal C}\times
{\cal C}$ with
\begin{equation}\label{def_dirichletform}
B^\theta_{\lambda ,\delta }({\bfm \varphi },{\bfm \psi })=\lambda
({\bfm \varphi},{\bfm \psi })_2 +(1/2)([{\bfm a}+{\bfm H}] D_x{\bfm
\varphi }, D_x{\bfm \psi })_2 -\theta( D_t{\bfm \varphi },{\bfm \psi
})_2+(\delta/2)(D_t{\bfm \varphi },D_t{\bfm \psi })_2.
\end{equation}
In what follows, the parameter $\theta$ (resp. $\delta$) is omitted
each time that it is equal to $1$ (resp. $0$). So the forms
$B^1_{\lambda,\delta}$, $B^\theta_{\lambda,0}$ and $B^1_{\lambda,0}$
are respectively simply denoted by $B_{\lambda,\delta}$,
$B^\theta_{\lambda}$ and $B_{\lambda}$.

\begin{proposition}\label{propdelta}
Suppose that ${\bfm h}\in L^2(\Omega )\cap {\rm Dom}(D_t)$ and
${\bfm d}\in{\cal H}$. Then, for any $\theta\in\{0;1\}$ and $\lambda
>0$, there exists a unique solution ${\bfm u}_\lambda \in \F$ of the
equation $\lambda {\bfm u}_\lambda-{\bfm L}{\bfm u}_\lambda- \theta
D_t{\bfm u}_\lambda={\bfm h}+{\bfm d}$,
%\begin{equation}\label{eq_res}
%\end{equation}
in the sense that $\forall {\bfm \varphi}\in \F$,
$B^\theta_\lambda({\bfm u}_\lambda ,{\bfm \varphi})=({\bfm h},{\bfm
\varphi})_2+\langle{\bfm d},{\bfm \varphi}\rangle_{-1,1}$. Moreover,
$D_t{\bfm u}_\lambda\in\H$ and
\begin{subequations}
\begin{equation}\label{estimates1}
  \lambda|{\bfm
u}_\lambda|_2^2+m\|{\bfm u}_\lambda\|_1^2\leq |{\bfm
h}|_2^2/\lambda+\|{\bfm d}|_{-1}^2/m,
\end{equation}
\begin{equation}\label{estimates2}
  \lambda|D_t{\bfm u}_\lambda|_2^2+m\|D_t{\bfm u}_\lambda\|_1^2\leq
|D_t{\bfm h}|_2^2/\lambda+2\|{\bfm
d}\|_T^2/m+2(C^a_2+C^H_2)^2\big(|{\bfm h}|^2/\lambda+\|{\bfm
d}\|_{-1}^2/m\big)/m^2.
\end{equation}
\end{subequations}
In the case ${\bfm d}\in L^2(\Omega)$, ${\bfm u}_\lambda\in {\rm
Dom}({\bfm L})$.\\ Finally, ${\bfm u}_\lambda$ is the strong limit
in $\H$ as $\delta$ goes to $0$ of the sequence $({\bfm
u}_{\lambda,\delta})_{\lambda,\delta}$, where ${\bfm
u}_{\lambda,\delta}$ is the unique solution of the equation:
$\forall {\bfm \varphi}\in \F$, $B^\theta_{\lambda,\delta}({\bfm
u}_{\lambda,\delta} ,{\bfm \varphi})=({\bfm h},{\bfm
\varphi})_2+\langle{\bfm d},{\bfm \varphi}\rangle_{-1,1}$, and the
family $(D_t{\bfm u}_{\lambda,\delta})_\delta$ is bounded in
$L^2(\Omega)$.
\end{proposition}

Before proving this result, we first investigate the case of time
independent coefficients. On the first side, this is a good starting
point for understanding the proof in the time dependent case and
this will bring out the difficulties arising with the time
dependency. On the other side, this result is needed in the last
section of this paper in order to prove the tightness of the process
$X$.

\begin{proposition}\label{existnotime}
Suppose that ${\bfm h}\in L^2(\Omega )$ Then, for any $\lambda
>0$, there exists a unique solution ${\bfm w}_\lambda \in \H\cap
{\rm Dom}({\bfm S})$ of the equation
\begin{equation}
\lambda {\bfm w}_\lambda-{\bfm S}{\bfm w}_\lambda={\bfm h}.
\end{equation}
\end{proposition}

\vspace{2mm} \noindent {\bf Proof :} The main tool of this proof is
the Lax-Milgram theorem. Let $\lambda >0$ be fixed. For any ${\bfm
\varphi },{\bfm \psi }\in {\cal C}$, consider the bilinear form on
${\cal C}\times {\cal C}$ defined by
$$D_\lambda ({\bfm \varphi },{\bfm \psi })=\lambda ({\bfm \varphi },{\bfm \psi })_2-({\bfm \varphi },{\bfm S}{\bfm \psi })_2.$$
Thanks to Assumption \ref{hypcontrol}, this form is clearly coercive
and continuous on ${\cal C}\times {\cal C}$ so that it can be
extended to the whole space $\H\times\H$.
%$$\forall {\bfm \varphi },{\bfm \psi }\in {\cal C},\quad D_\lambda ({\bfm \varphi },{\bfm \psi })\leq \max(\lambda ,M)\|{\bfm \varphi }\|_1\|{\bfm \psi }\|_1$$ and also
%$$\forall {\bfm \varphi }\in {\cal C},\quad \min(\lambda ,m)\|{\bfm \varphi }\|^2_1\leq D_\lambda ({\bfm \varphi },{\bfm \varphi }).$$
 The extension is also coercive and continuous. Now, the application ${\bfm \varphi }\mapsto ({\bfm h },{\bfm \varphi })_2 $ is
 obviously continuous on $\H$ so that the Lax-Milgram theorem
 applies. It allows to construct a strongly continuous resolvent associated to $\lambda-{\bfm
 S}$ by way of classical tools (see \cite[Ch. 1, Sect. 3]{fukushima} or \cite[Ch. 1, Sect. 2]{ma} for further details).\qed

\vspace{2mm} \noindent {\bf Proof of the Proposition
\ref{propdelta}:} Since the case $\theta =0$ and $\theta=1$ are
quite similar, we only give the proof for $\theta=1$. The existence
of a solution relies on the Lax-Milgram theorem again. However, the
considered bilinear form \eqref{def_dirichletform} with $\delta=0$
is not coercive on $\F$ because of the time differential term
$(D_t{\bfm \varphi },{\bfm \psi })$. The strategy consists in making
it coercive by adding a term $(\delta/2)(D_t{\bfm \varphi },D_t{\bfm
\psi })$ ($\delta>0$) and then letting $\delta$ go to $0$. Notice
that for ${\bfm \varphi },{\bfm \psi }\in {\cal C}$, we have
$$\big([\lambda -{\bfm L}- D_t-( \delta /2)D_t^2]{\bfm
\varphi},{\bfm \psi }\big)_2=B_{\lambda ,\delta } ({\bfm \varphi
},{\bfm \psi }).$$ The continuity of $B_{\lambda ,\delta }$ on
${\cal C}\times {\cal C}\subset \F\times \F$ follows from
(\ref{inega}) and (\ref{inegh}). As a result of the
time-independence of ${\bfm V}$, for any ${\bfm \varphi }\in {\cal
C}$, we have $({\bfm \varphi },D_t{\bfm \varphi })_2=0 $. As a
consequence, for any ${\bfm \varphi }\in {\cal C}$,
\begin{equation}\label{coer}
\min(\lambda ,\delta /2,m)\varepsilon({\bfm \varphi }, {\bfm \varphi
})\leq B_{\lambda ,\delta }({\bfm \varphi },{\bfm \varphi }).
\end{equation}
Hence $B_{\lambda ,\delta }$ defines a continuous coercive bilinear
form on $\F\times \F$. The Lax-Milgram theorem applies and provides
us with a solution ${\bfm u}_{\lambda ,\delta }$ of the equation
\begin{equation}\label{eqlamdel}
\forall {\bfm \varphi }\in {\cal C},\quad B_{\lambda ,\delta }({\bfm
u}_{\lambda,\delta } ,{\bfm \varphi })=({\bfm h}, {\bfm \varphi
})_2+\langle{\bfm d},{\bfm \varphi}\rangle_{-1,1}.
\end{equation}
In particular, choosing ${\bfm \varphi}={\bfm u}_{\lambda,\delta}$
in \eqref{eqlamdel}, we get the bound
\begin{equation}\label{bound_h}
 \lambda|{\bfm
u}_{\lambda,\delta }|_2^2+m\|{\bfm u}_{\lambda,\delta
}\|_1^2+\delta|D_t{\bfm u}_{\lambda,\delta }|_2^2\leq  |{\bfm
h}|_2^2/\lambda+\|{\bfm d}\|_{-1}^2/m.
\end{equation}
Let us now to pass to the limit as $\delta $ goes to $0$ to obtain a
solution ${\bfm u}_\lambda\in\F $ of the equation
\begin{equation}\label{eqlam}
\forall {\bfm \varphi }\in {\cal C},\quad B_\lambda ({\bfm
u}_\lambda ,{\bfm \varphi })=({\bfm h}, {\bfm \varphi
})_2+\langle{\bfm d},{\bfm \varphi}\rangle_{-1,1} .
\end{equation}
We are faced with the problem of controlling $D_t{\bfm
u}_{\lambda,\delta }$ as $\delta$ goes to $0$. The idea lies in
differentiating equation (\ref{eqlamdel}) with respect to the time
variable in order to establish an equation satisfied by $D_t{\bfm
u}_{\lambda,\delta }$, from which estimates will be derived. So, we
define for each fixed $\lambda ,\delta>0$, ${\bfm
v}_{s}=\Lambda_s{\bfm u}_{\lambda ,\delta }$ (the parameters
$\lambda,\delta$ of ${\bfm v}_s$ are temporarily omitted in order to
simplify the notations) and we easily check that ${\bfm v}_{s}$
solves the following equation
\begin{equation}\label{resdtt}
\forall {\bfm \varphi }\in \F,\,B_{\lambda ,\delta }({\bfm
v}_{s},{\bfm \varphi })={\bfm F}_s({\bfm \varphi }),
\end{equation}
where ${\bfm F}_s$ is a continuous linear form on $\F$ defined,
$\forall {\bfm \varphi }\in \F$, by
\begin{equation}
{\bfm F}_s({\bfm \varphi }) = (\Lambda_s{\bfm h} , {\bfm \varphi
})_2-\langle{\bfm d},\Lambda_{-s}{\bfm \varphi}\rangle_{-1,1}
-[\Lambda_s{\bfm T}_a+\Lambda_s{\bfm T}_H]\big((-\widetilde{{\bfm
S}})^{1/2}T_{s,0}{\bfm u}_{\lambda ,\delta },(-\widetilde{{\bfm
S}})^{1/2}{\bfm \varphi }\big).
\end{equation}
From Assumption \ref{hypcontrol}, it is readily seen that $${\bfm
F}_s({\bfm \varphi })\leq |D_t{\bfm h}|_2|{\bfm \varphi}|_2+\|{\bfm
d}\|_T\|{\bfm \varphi }\|_{1}+(C^a_2+C_2^H)\|{\bfm u}_{\lambda
,\delta }\|_1\|{\bfm \varphi }\|_{1},$$ for any $ s\in\R^*$.
Therefore
\begin{equation}\label{bornebvs}
B_{\lambda ,\delta }({\bfm v}_s,{\bfm v}_s) = {\bfm F}_s({\bfm
v}_s)\leq |D_t{\bfm h}|_2|{\bfm v}_s|_2+\|{\bfm d}\|_T\|{\bfm
\varphi }\|_{1}+(C^a_2 +C_2^H)\|{\bfm u}_{\lambda ,\delta
}\|_1\|{\bfm v}_s\|_{1}.
 \end{equation}
Using estimate \eqref{bound_h} in \eqref{bornebvs}, we have
\begin{equation}\label{bound_dt} \lambda|{\bfm
v}_s|_2^2+m\|{\bfm v}_s\|_{1}^2+\delta|D_t{\bfm v}_s|_{2}^2\leq
|D_t{\bfm h}|_2^2/\lambda+2\|{\bfm
d}\|_T^2/m+2(C^a_2+C^H_2)^2\big(|{\bfm h}|^2/\lambda+\|{\bfm
d}\|_{-1}^2/m\big)/m^2.
\end{equation}
So, the family $({\bfm v}_s)_{s\in \R^*}$ is bounded in $\F$. Even
if it means extracting a subsequence (still denoted by $({\bfm
v}_s)_{s\in \R^*}$), $({\bfm v}_s)_{s\in \R^*}$ converges weakly in
$\F$ towards some function ${\bfm v}_0\in\F$ as $s$ goes to $0$. On
the other hand, since ${\bfm u}_{\lambda,\delta}\in\F\subset{\rm
Dom}(D_t)$, $({\bfm v}_s)_{s\in \R^*}$ also converges strongly in
$L^2(\Omega )$ towards $D_t{\bfm u}_{\lambda ,\delta }$, so that
$D_t{\bfm u}_{\lambda ,\delta }\in \F$ and satisfies bound
\eqref{bound_dt} instead of ${\bfm v}_s$. In particular, $(D_t{\bfm
u}_{\lambda,\delta})_{\delta>0}$ is bounded in $\H$ independently of
$\delta>0$ and so is $({\bfm u}_{\lambda,\delta})_{\delta>0}$ in
$\F$. %As a by-product, if ${\bfm d}=0$, note that $D_t{\bfm
%u}_{\lambda,\delta}$ then solves the equation
%\begin{equation}\label{eq_dtuld}
%\forall {\bfm \varphi}\in \F,\quad B_{\lambda,\delta}(D_t{\bfm
%u}_{\lambda,\delta},{\bfm \varphi})=(D_t{\bfm h},{\bfm
%\varphi})_2+[\partial_t{\bfm T}_a+\partial_t{\bfm
%T}_H]\big((-\widetilde{{\bfm S}})^{1/2}{\bfm
%u}_{\lambda,\delta},(-\widetilde{{\bfm S}})^{1/2}{\bfm
%\varphi}\big).
%\end{equation}
Thereby, there exists a subsequence $({\bfm u}_{\lambda,\delta
},D_t{\bfm u}_{\lambda,\delta })_{\delta >0}\subset \F\times\H$,
still indexed with $\delta>0$, that converges weakly in $\F\times\H$
towards $({\bfm u}_\lambda,D_t{\bfm u}_\lambda)\in \F\times\H$  as
$\delta \rightarrow 0$. In particular, $\delta D_t{\bfm u}_{\lambda
,\delta }\rightarrow 0$ in $L^2(\Omega)$ as $\delta$ goes to $0$. So
we are in position to pass to the limit as $\delta$ goes to $0$ in
\eqref{eqlamdel}. Obviously, ${\bfm u}_\lambda$ is a solution of
\eqref{eqlam}. Uniqueness of the weak limit raises no particular
difficulty since two weak limits ${\bfm u}_\lambda $ and ${\bfm
w}_\lambda $ satisfy $\forall {\bfm \varphi }\in \F$, $B_\lambda
({\bfm u}_\lambda-{\bfm w}_\lambda,{\bfm \varphi })=0$. It just
remains to choose ${\bfm \varphi }= {\bfm u}_\lambda-{\bfm
w}_\lambda$. \eqref{estimates1} and \eqref{estimates2} respectively
result from \eqref{bound_h} and \eqref{bound_dt}. If ${\bfm d}\in
L^2(\Omega)$, note that ${\bfm u}_\lambda\in\F\subset\H$ and ${\cal
B}_\lambda
({\bfm u}_\lambda,\cdot)=({\bfm h}+{\bfm d}+D_t{\bfm u}_\lambda,\cdot)_2$ is $L^2$-continuous so that ${\bfm u}_\lambda\in {\rm Dom}({\bfm L})$. \\
Let us now investigate the strong convergence in $\F$ of $({\bfm
u}_{\lambda,\delta})_{\lambda,\delta}$ towards ${\bfm u}_{\lambda}$
as $\delta$ goes to $0$. Let us make the difference between
\eqref{eqlamdel} and \eqref{eqlam} and choose ${\bfm \varphi}={\bfm
u}_{\lambda,\delta}-{\bfm u}_{\lambda}$, this yields
$$B_{\lambda,\delta}({\bfm
u}_{\lambda,\delta}-{\bfm u}_{\lambda},{\bfm
u}_{\lambda,\delta}-{\bfm u}_{\lambda}) =(\delta/2)(D_t{\bfm
u}_{\lambda},D_t{\bfm u}_{\lambda}-D_t{\bfm
u}_{\lambda,\delta})_2,$$ and this latter quantity converges to $0$
as $\delta$ goes to $0$ because of the boundedness of the family
$(|D_t{\bfm u}_{\lambda,\delta}|_2)_{\lambda,\delta}$. \eqref{coer}
allows to conclude. %Note that $D_t{\bfm u}_\lambda$ satisfies
%equation
%\begin{equation}\label{eq_derive}
%\begin{split}
%\lambda  |D_t{\bfm u}_{\lambda}|_2^2+{\bfm
%T}_a\big((-\widetilde{{\bfm S}})^{1/2}D_t{\bfm
%u}_{\lambda},(-\widetilde{{\bfm S}})^{1/2}D_t{\bfm
%u}_{\lambda}\big)= & (D_t{\bfm u}_{\lambda},D_t{\bfm h})_2\\
%+ & [\partial_t{\bfm T}_a+\partial_t{\bfm
%T}_H]\big((-\widetilde{{\bfm S}})^{1/2}{\bfm
%u}_{\lambda,\delta},(-\widetilde{{\bfm S}})^{1/2}D_t{\bfm
%u}_{\lambda}\big).
%\end{split}
%\end{equation}
\qed

\subsection{Control of the solution}\label{sec_control}
%%%%%%%%%%%%%%%%%%%%%%%%%%%%%%%%%%%%%%%%%%%%%%%%%%%%%%%%%%%%%%%%%%%%%%%%%%%%%%%%%%%%%%%%%%%%%%%%%%%%%%%%

Our goal is now to determine the asymptotic behaviour, as $\lambda$
goes to $0$, of the solution ${\bfm u}_\lambda ^i$ of the equation
(in the sense of Proposition \ref{propdelta})
\begin{equation}\label{eqb}
\lambda {\bfm u}_\lambda ^i-{\bfm L}{\bfm u}_\lambda ^i-D_t{\bfm
u}_\lambda ^i={\bfm b}_i.
\end{equation}
More precisely, we aim at proving that $\lambda |{\bfm u}_\lambda
^i|_2^2\rightarrow 0$ and that $(\nabla ^\sigma {\bfm u}_\lambda
^i)_\lambda $ converges in $(L^2(\Omega ))^d$ as $\lambda $ goes to
$0$. Our strategy consists in showing that the operator $\lambda
-{\bfm L} -D_t$  is just a perturbation of the operator $\lambda
-\widetilde{{\bfm S}} -D_t$,  so that the study can be reduced to
studying the solution of the equation
$$\lambda {\bfm v}_\lambda -\widetilde{{\bfm S}}{\bfm v}_\lambda -D_t{\bfm v}_\lambda={\bfm b}_\lambda ,$$ where ${\bfm b}_\lambda $
will be defined thereafter but possesses a strong limit in
$\H_{-1}$. This latter equation is more convenient to study because
the operators $\widetilde{{\bfm S}}$ and $D_t$ can be viewed through
the same spectral decomposition. Thus, the purpose of this section
is to prove the following Proposition

\begin{proposition}\label{propcon}
Let $({\bfm b}_\lambda )_{\lambda>0}$ be a family of functions in
$\H_{-1}\cap L^2(\Omega )$ which is strongly convergent in $\H_{-1}$
to ${\bfm b}_0$. Suppose that there exists a constant $C$ (which
does not depend on $\lambda $) such that $\forall s>0$ and $\forall
{\bfm \varphi }\in {\cal C}$,
$$({\bfm b}_\lambda ,\Lambda_s{\bfm \varphi })_2\leq
C \| {\bfm \varphi }\|_1.$$ Then the solution ${\bfm u}_\lambda \in
\F$ of the equation $\lambda {\bfm u}_\lambda -{\bfm L}{\bfm
u}_\lambda -D_t{\bfm u}_\lambda ={\bfm b}_\lambda $ (in the sense of
Proposition \ref{propdelta}) satisfies:

$\bullet$ there exists ${\bfm \eta }\in \D$ such that
$(-\widetilde{{\bfm S}})^{1/2}{\bfm u}_\lambda \rightarrow {\bfm
\eta }$ as $\lambda$ goes to $0$ in $\D$,

$\bullet$  $\lambda |{\bfm u}_\lambda|_2^2\rightarrow 0$ as
$\lambda$ goes to $0$.
\end{proposition}

As for the existence of the solution, let us first investigate the
time independent case by way of introduction.

\begin{proposition}\label{controlnodt}
Let ${\bfm h}$ be in $ \H_{-1}\cap L^2(\Omega )$. For any $\lambda
>0$, let ${\bfm w}_\lambda$ be defined as the unique solution in
$\H$ of the equation
$$\lambda   {\bfm w}_\lambda -{\bfm S}{\bfm w}_\lambda ={\bfm h}$$
Then $\lambda |{\bfm w}_\lambda|_2^2\rightarrow 0$ and there exists
${\bfm \zeta }\in (L^2(\Omega ))^d$ such that $|\nabla ^\sigma {\bfm
w}_\lambda-{\bfm \zeta }|_2\rightarrow 0$ as $\lambda$ goes to $0$.
\end{proposition}
\noindent {\bf Proof :} Keeping the notations of Proposition
\ref{existnotime}, ${\bfm w}_\lambda$ solves the equation: $\forall
{\bfm \varphi}\in \H$, $D_\lambda({\bfm w}_\lambda,{\bfm
\varphi})=({\bfm h},{\bfm \varphi})_2$. Choosing ${\bfm
\varphi}={\bfm w}_\lambda$ and using ${\bfm h}\in\H_{-1}$, we have
$\lambda |{\bfm w}_\lambda|_2^2+m\|{\bfm w}_\lambda\|_1^2\leq
\|{\bfm h}\|_{-1}^2/m$. Thus, even if it means extracting a
subsequence, we can find ${\bfm g}\in L^2(\Omega )$ such that
$((-\widetilde{{\bfm S}})^{1/2}{\bfm w}_\lambda)_\lambda $ converges
weakly in $L^2(\Omega )$ towards ${\bfm g}$ as $\lambda$ tends to
$0$. Moreover $(\lambda {\bfm w}_\lambda)_\lambda$ clearly converges
to $0$ in $L^2(\Omega )$. For any ${\bfm \varphi }\in \H$, passing
to the limit as $\lambda$ goes to zero in the expression
\begin{equation}
\lambda ({\bfm w}_\lambda,{\bfm \varphi })_2 +{\bfm
T}_a\big((-\widetilde{{\bfm S}})^{1/2} {\bfm
w}_\lambda,(-\widetilde{{\bfm S}}\big)^{1/2}{\bfm \varphi })_2
=D_\lambda ({\bfm w}_\lambda,{\bfm \varphi })=({\bfm h},{\bfm
\varphi })_2 ,
\end{equation}
we obtain ${\bfm T}_a\big( {\bfm g},(-\widetilde{{\bfm
S}})^{1/2}{\bfm \varphi }\big)_2=({\bfm h},{\bfm \varphi })_2 $.
%\begin{equation}\label{scalare}
%\end{equation}
Making the difference between the last two equalities, subtracting
${\bfm T}_a\big((-\widetilde{{\bfm S}})^{1/2}{\bfm w}_\lambda-{\bfm
g},{\bfm g}\big) $ and then choosing $(-\widetilde{{\bfm
S}})^{1/2}{\bfm \varphi}=(-\widetilde{{\bfm S}})^{1/2}{\bfm
w}_\lambda-{\bfm g}$, we obtain
\begin{eqnarray*}
\lambda|{\bfm w}_\lambda|_2^2 +{\bfm T}_a\big((-\widetilde{{\bfm
S}})^{1/2}{\bfm w}_\lambda-{\bfm g},(-\widetilde{{\bfm
S}})^{1/2}{\bfm w}_\lambda-{\bfm g}\big)=-{\bfm
T}_a\big((-\widetilde{{\bfm S}})^{1/2}{\bfm w}_\lambda-{\bfm
g},{\bfm g}\big).
\end{eqnarray*}
Due to the weak convergence of $((-\widetilde{{\bfm S}})^{1/2}{\bfm
w}_\lambda)_\lambda$ to ${\bfm g}$ in $\D$, the right-hand side
converges to $0$ as $\lambda$ goes to $0$. So does the left-hand
side. Since ${\bfm T}_a$ defines an inner product on $\D$ equivalent
to the canonical one (Assumption \ref{hypcontrol}), this completes
the proof of the strong convergence up to a subsequence. Uniqueness
of the weak limit is clear since two weak limits ${\bfm g}$ and
${\bfm g}'\in\D$ satisfy: $\forall{\bfm \varphi}\in {\cal C}$,
${\bfm T}_a({\bfm g}-{\bfm g}',(-\widetilde{{\bfm S}})^{1/2}{\bfm
\varphi})=0$. Finally, since the convergence in $\D$ of
$((-\widetilde{{\bfm S}})^{1/2}{\bfm w}_\lambda)_\lambda$ is
equivalent to the convergence of $(\nabla^\sigma{\bfm
w}_\lambda)_\lambda$ in $(L^2(\Omega))^d $, we complete the
proof.\qed

\begin{proposition}\label{convsdt}
Let $({\bfm b}_\lambda )_{\lambda >0}$ be a family of functions in
$\H_{-1}$ that is strongly convergent to ${\bfm b}_0$ in $\H_{-1}$.
Let $({\bfm v}_\lambda )_{\lambda >0}$ be a family of functions in
$\F$  that solves the equation (for any $ \lambda >0$) $\lambda
{\bfm v}_\lambda -\widetilde{{\bfm S}}{\bfm v}_\lambda- D_t {\bfm
v}_\lambda ={\bfm b}_\lambda$ in the following sense,
\begin{equation}\label{resv}
\forall {\bfm \varphi}\in \F,\quad \lambda({\bfm v}_\lambda,{\bfm
\varphi})_2+\langle{\bfm v}_\lambda,{\bfm
\varphi}\rangle_1-(D_t{\bfm v}_\lambda,{\bfm \varphi})_2=({\bfm
b}_\lambda,{\bfm \varphi})_2.
\end{equation}
Then there exists ${\bfm \eta }\in \D$ such that $\lambda |{\bfm
v}_\lambda |_2^2\rightarrow 0$ and $|(-\widetilde{{\bfm S}})^{1/2}
{\bfm v}_\lambda-{\bfm \eta }|_{2}\rightarrow 0$ as $\lambda$ goes
to $0$.
\end{proposition}
\noindent {\bf Proof:} From Lemma \ref{lemma1} and Lemma
\ref{lemma2} below, we can assume that, for any $\lambda>0$, ${\bfm
b}_\lambda\in L^2(\Omega)\cap {\rm Dom}(D_t)\cap\H_{-1}$ and
converges to ${\bfm b}_0\in\H_{-1}$. Then ${\bfm v}_\lambda\in{\rm
Dom}(\widetilde{{\bfm S}})$ (see Proposition \ref{propdelta}).
Remind that $-\widetilde{{\bfm S}}=\int_{\reel _+\times \reel
}x\,E(dx,dy)$ and $-D_t= \int_{\reel _+\times \reel }iy\,E(dx,dy)$.
Choosing ${\bfm \varphi}={\bfm v}_\lambda$ in \eqref{resv}, we have
\begin{equation}\label{env}
\lambda |{\bfm v}_\lambda |_2^2+\|{\bfm v}_\lambda\|_{1}^2 =({\bfm
b}_\lambda, {\bfm v}_\lambda )_2 \leq C \|{\bfm v}_\lambda\|_{1}\leq
C^2,
\end{equation}
where $C=\sup_{\lambda >0}\|{\bfm b}_\lambda \|_{-1}$. Thus we can
find  ${\bfm h}\in \D$ and a subsequence, still denoted by $( {\bfm
v}_\lambda)_\lambda $, such that $\left((-\widetilde{{\bfm
S}})^{1/2} {\bfm v}_\lambda\right)_\lambda $ converges weakly in
$L^2(\Omega )$ to ${\bfm h}$.\\
Now we claim $\sup_{\lambda >0}\|\lambda {\bfm v}_\lambda
\|_{{-1}}<\infty $ and $ \sup_{\lambda
>0}\|D_t {\bfm v}_\lambda
\|_{{-1}}<\infty $.
\begin{eqnarray*}
|(\lambda {\bfm v}_\lambda ,{\bfm \varphi })_2 | & = &
\big|\int_{\reel _+\times \reel } \lambda (\lambda
+x+iy)^{-1}\,dE_{{\bfm b}_\lambda ,{\bfm \varphi }}\big|\\ & \leq &
\Big(\int_{\reel _+\times \reel }\frac{ \lambda ^2}{ x[(\lambda
+x)^2+y^2]}\,dE_{{\bfm b}_\lambda,{\bfm
b}_\lambda}\Big)^{1/2}\Big(\int_{
\reel _+\times \reel } x\,dE_{{\bfm \varphi },{\bfm \varphi }}\Big)^{1/2}\\
& \leq  & \sup_{\lambda >0}\Big( \int_{\reel _+\times \reel }
x^{-1}\,dE_
{{\bfm b}_\lambda,{\bfm b}_\lambda}\Big)^{1/2}\|{\bfm \varphi }\|_1\\
& =& \sup_{\lambda >0}\|{\bfm b}_\lambda\|_{{-1}}\|{\bfm \varphi
}\|_1.
%\\ & \leq & \overline{C}\|{\bfm \varphi }\|_1
\end{eqnarray*}
%where $\overline{C}$ does not depend on $\lambda $.
Since $D_t{\bfm v}_\lambda =\lambda {\bfm v}_\lambda
-\widetilde{{\bfm S}}{\bfm v}_\lambda -{\bfm b}_\lambda  $ and
$\|\widetilde{{\bfm S}}{\bfm v}_\lambda\|_{-1}\leq \|{\bfm
v}_\lambda\|_{1}$, $D_t{\bfm v}_\lambda \in \H_{-1}$ and
$\sup_{\lambda>0}\|D_t{\bfm v}_\lambda\|_{-1}<~\infty$. Then there
exists a bounded family $({\bfm F}_\lambda )_{\lambda \geq 0}$ of
continuous linear forms on $\D\subset L^2(\Omega )$ such that
$\forall \lambda
>0$, $\forall {\bfm \varphi }\in {\cal C},\quad {\bfm F}_\lambda
((-\widetilde{{\bfm S}})^{1/2}{\bfm \varphi })=(D_t{\bfm v}_\lambda,
{\bfm \varphi })_2$. Moreover, from \eqref{env}, $(\lambda {\bfm
v}_\lambda)_\lambda$ converges to $0$ in $L^2(\Omega)$ so that,
$\forall {\bfm \varphi }\in {\cal C}$
\begin{eqnarray*}
{\bfm F}_\lambda ((-\widetilde{{\bfm S}})^{1/2}{\bfm \varphi }) & =
& (\lambda {\bfm v}_\lambda,{\bfm \varphi })_2+((-\widetilde{{\bfm
S}})^{1/2}{\bfm v}_\lambda,(-\widetilde{{\bfm S}})^{1/2}{\bfm
\varphi })_2-\langle{\bfm b}_\lambda,{\bfm \varphi}\rangle_{-1,1}\\
& \rightarrow & ({\bfm h},(-\widetilde{{\bfm S}})^{1/2}{\bfm \varphi
})_2-\langle{\bfm b}_0,{\bfm \varphi}\rangle_{-1,1}
\end{eqnarray*}
 as $\lambda$ goes to $0$. Hence, $({\bfm F}_\lambda )_{\lambda \geq
 0}$ is weakly convergent in $\D^*$ (topological dual of $\D$) to a
 limit denoted by ${\bfm F}_0$.\\
We now aim at proving $ {\bfm F}_0 ( {\bfm h} )=0 $. Using the
antisymmetry of the operator $D_t$
$${\bfm F}_\lambda ((-\widetilde{{\bfm S}})^{1/2}{\bfm v}_\mu ) = (D_t{\bfm v}_\lambda ,{\bfm v}_\mu)_2 = -( D_t{\bfm v}_\mu {\bfm v}_\lambda )_2
 = -{\bfm F}_\mu  ((-\widetilde{{\bfm S}})^{1/2}{\bfm v}_\lambda  ),$$
we pass to the limit as $\lambda $ goes to $0$ and obtain ${\bfm
F}_0((-\widetilde{{\bfm S}})^{1/2}{\bfm v}_\mu )=-{\bfm F}_\mu (
{\bfm h} )$. It just remains to pass to the limit as $\mu$ goes to
$0$, it yields ${\bfm F}_0({\bfm h})=-{\bfm F}_0({\bfm h })=0$.\\
Let us investigate now the limit equation, which connects ${\bfm
F}_0$, ${\bfm h}$ and ${\bfm b}_0$. First remind of (\ref{env}),
which states $\lambda |{\bfm v}_\lambda |_2^2\leq C^2$ and as a
consequence $\lambda {\bfm v}_\lambda\rightarrow 0$ as $\lambda$
goes to $0$. Then, we are in a position to pass to the limit as
$\lambda$ tends to $0$ in (\ref{resv}), and this yields, for any
${\bfm \varphi }\in \F$,
\begin{equation}\label{limitequ}
({\bfm h},(-\widetilde{{\bfm S}})^{1/2}{\bfm \varphi })_2 -{\bfm
F}_0 ((-\widetilde{{\bfm S}})^{1/2} {\bfm \varphi }) =\langle{\bfm
b}_0,{\bfm \varphi }\rangle_{-1,1}.
\end{equation}
Let us now establish the uniqueness of the weak limit. Let ${\bfm h
}$ and ${\bfm h}'$ be two possible weak limits of two subsequences
of $({\bfm v}_\lambda)_{\lambda }$, and ${\bfm F}_0$,${\bfm F}'_0$
the corresponding linear forms defined as described above. Then
\eqref{limitequ} provides us with he following equality:
\begin{equation}\label{eqlim1}
\forall{\bfm \varphi }\in \F,\quad ({\bfm h}-{\bfm
h}',(-\widetilde{{\bfm S}})^{1/2}{\bfm \varphi }) =[{\bfm F}_0-{\bfm
F}_0']((-\widetilde{{\bfm S}})^{1/2}{\bfm \varphi }).
\end{equation}
Using the antisymmetry of the operator $D_t$ again, we obtain
$${\bfm F}_\lambda ((-\widetilde{{\bfm S}})^{1/2}{\bfm v}_\mu ) = (D_t{\bfm v}_\lambda ,{\bfm v}_\mu )_2 = -( D_t{\bfm v}_\mu ,{\bfm v}_\lambda )_2
 = -{\bfm F}_\mu ((-\widetilde{{\bfm S}})^{1/2}{\bfm v}_\lambda  ).$$
Let us first pass to the limit as $\lambda $ goes to $0$ along the
first subsequence, and then pass to the limit as $\mu $ goes to $0$
along the second subsequence, we obtain
\begin{equation*}\label{signf}
{\bfm F}_0 ({\bfm h}' )=-{\bfm F}'_0({\bfm h } ).
\end{equation*}
Now, it just remains to choose $(-\widetilde{{\bfm S}})^{1/2}{\bfm
\varphi }={\bfm h}- {\bfm h}'$ in (\ref{eqlim1}) and this yields
\begin{equation*}
 |{\bfm h}-{\bfm h'}|_2^2=-{\bfm F}_0({\bfm h'}) - {\bfm
F}'_0({\bfm h}) =0.
\end{equation*} Hence the
weak convergence holds for the whole family. Let us now tackle the
strong convergence of $({\bfm v}_\lambda)_\lambda$. Choosing ${\bfm
\varphi}={\bfm v}_\lambda$ in \eqref{limitequ}, using ${\bfm
F}_0({\bfm h })=0$ and passing to the limit a $\lambda$ goes to $0$,
this yields
\begin{equation}\label{convh}
 ( {\bfm h },{\bfm h } )_2 = \lim_{\lambda
\to 0} \langle {\bfm b }_0,{\bfm v}_\lambda
\rangle_{-1,1}=\lim_{\lambda \to 0} \langle {\bfm b }_\lambda,{\bfm
v}_\lambda \rangle_{-1,1} =\lim_{\lambda \to 0} \big[\lambda|{\bfm v
}_\lambda|_2^2 +\|{\bfm v}_\lambda\|_1^2 \big].
\end{equation}
In particular, $|{\bfm h }|_2=\lim_{\lambda \to 0}
|(-\widetilde{{\bfm S}})^{1/2}{\bfm v}_\lambda|_2$.  Thus, the
convergence of the norms implies the strong convergence of the
sequence $((-\widetilde{{\bfm S}})^{1/2}{\bfm v}_\lambda)_\lambda$
to ${\bfm h}$ in $L^2(\Omega)$. As a bypass, \eqref{convh} also
implies the convergence of $\big(\lambda |{\bfm
v}_\lambda|_2^2\big)_\lambda$ to $0$. \qed

\begin{lemma}\label{lemma1}
For each function ${\bfm b}\in \H_{-1}$, there exists a family
$({\bfm b}_\lambda)_\lambda$ of functions in $L^2(\Omega)\cap{\rm
Dom}(D_t)\cap \H_{-1}$ such that $\|{\bfm b}-{\bfm
b}_\lambda\|_{-1}$ converges to $0$ as $\lambda$ goes to $0$.
\end{lemma}

\noindent {\bf Proof:} Let us consider the solution ${\bfm
w}_\lambda\in \H$ of the equation $\lambda {\bfm
w}_\lambda-\widetilde{{\bfm S}}{\bfm w}_\lambda={\bfm b}$ (see
Proposition \ref{existnotime}). Then, for any ${\bfm \varphi}\in
{\cal C}$,
\begin{eqnarray*}
(\lambda {\bfm w}_\lambda,{\bfm \varphi})_2 & = &
\int_{\R^+\times\R}\lambda(\lambda+x)^{-1}\,dE_{{\bfm b},{\bfm
\varphi}}(dx,dy)\\ & \leq &
\Big(\int_{\R^+\times\R}\lambda^2x^{-1}(\lambda+x)^{-2}\,dE_{{\bfm
b},{\bfm b}}(dx,dy)\Big)^{1/2}\|{\bfm \varphi}\|_1.
\end{eqnarray*}
Since ${\bfm b}\in\H_{-1}$, we have
$\int_{\R^+\times\R}x^{-1}\,dE_{{\bfm b},{\bfm b}}(dx,dy)<\infty$.
Thus the Lebesgue theorem ensures that the above integral converges
to $0$ as $\lambda$ goes to $0$. Hence, $\|\lambda {\bfm
w}_\lambda\|_{-1}$ converges to $0$ as $\lambda$ goes to $0$. We can
now choose a family $({\bfm \varphi}_\lambda)_\lambda$ in ${\cal C}$
such that $\|{\bfm w}_\lambda-{\bfm \varphi}_\lambda\|_1\rightarrow
0$ as $\lambda$ goes to $0$. Finally,
$$ \|{\bfm
b}-\widetilde{{\bfm S}}{\bfm \varphi}_\lambda\|_{-1}\leq  \|{\bfm
b}-\widetilde{{\bfm S}}{\bfm w}_\lambda\|_{-1}+ \|\widetilde{{\bfm
S}}{\bfm w}_\lambda-\widetilde{{\bfm S}}{\bfm
\varphi}_\lambda\|_{-1}\leq \|\lambda {\bfm
w}_\lambda\|_{-1}+\|{\bfm w}_\lambda-{\bfm \varphi}_\lambda\|_1$$
also converges to $0$ as $\lambda$ tends to $0$ and, clearly,
$\widetilde{{\bfm S}}{\bfm \varphi}_\lambda\in L^2(\Omega)\cap {\rm
Dom}(D_t)$.\qed

\begin{lemma}\label{lemma2}
Let $({\bfm b}_\lambda)_\lambda$ and $({\bfm b}'_\lambda)_\lambda$
be two families in $\H_{-1}$ such that $\|{\bfm b}_\lambda-{\bfm
b}'_\lambda\|_{-1}\rightarrow 0$ as $\lambda$ goes to $0$. Let
$({\bfm v}_\lambda)_\lambda$ and $({\bfm v}'_\lambda)_\lambda$ two
families in $\F$ solving equation \eqref{resv} with respectively
${\bfm b}_\lambda$ and ${\bfm b}'_\lambda$ as right-hand side. Then
$\lambda|{\bfm v}_\lambda-{\bfm v}_\lambda'|_2^2+\|{\bfm
v}_\lambda-{\bfm v}_\lambda'\|_1^2\rightarrow 0$ as $\lambda$ goes
to $0$.
\end{lemma}

\noindent {\bf Proof:} Making the difference between the two
equations corresponding to ${\bfm v}_\lambda$ and ${\bfm
v}'_\lambda$, this yields for any ${\bfm \varphi}\in \F$,
$$\lambda({\bfm v}_\lambda-{\bfm v}_\lambda',{\bfm
\varphi})_2+\langle {\bfm v}_\lambda-{\bfm v}_\lambda',{\bfm
\varphi}\rangle_1-(D_t{\bfm v}_\lambda-D_t{\bfm v}_\lambda',{\bfm
\varphi})_2=\langle {\bfm b}_\lambda-{\bfm b}_\lambda',{\bfm
\varphi}\rangle_{-1,1}.$$ Choosing ${\bfm \varphi}={\bfm
v}_\lambda-{\bfm v}_\lambda'$, we easily deduce $\lambda|{\bfm
v}_\lambda-{\bfm v}_\lambda'|_2^2+\|{\bfm v}_\lambda-{\bfm
v}_\lambda'\|_1\leq \|{\bfm b}_\lambda-{\bfm b}_\lambda'\|_{-1}$.
The result follows.\qed

 Let us now investigate
the general case, that means that we aim at replacing
$\widetilde{{\bfm S}}$ by ${\bfm L}$ in Proposition \ref{convsdt}.
We first set out the main ideas of the proof. Let us formally write
\begin{eqnarray*}
\lambda -{\bfm L}-D_t & = & \lambda -\widetilde{{\bfm S}}-D_t -({\bfm L}-\widetilde{{\bfm S}})\\
 & = & \big({\rm I}-\big[{\bfm L}-\widetilde{{\bfm S}}\big](\lambda -\widetilde{{\bfm S}}-D_t)^{-1}\big)(\lambda -\widetilde{{\bfm S}}-D_t)
\end{eqnarray*}
If we can prove that $\big[{\bfm L}-\widetilde{{\bfm
S}}\big](\lambda -\widetilde{{\bfm S}}-D_t)^{-1}$ defines a strictly
contractive operator, then we will be in position to inverse it. It
turns out that it is actually bounded but not strictly contractive.
To overcome this difficulty, we introduce a small parameter $\delta
$ to make the operator $\delta \big[{\bfm L}-\widetilde{{\bfm
S}}\big](\lambda -\widetilde{{\bfm S}}-D_t)^{-1}$  strictly
contractive. Then, an iteration procedure proves that $\delta$ can
be chosen equal to $1$.

\begin{proposition}\label{rec}
Let $({\bfm b}_\lambda )_{\lambda>0}$ be a family of functions in
$\H_{-1}$ that  is strongly convergent in $\H_{-1}$ to some ${\bfm
b}_0\in\H_{-1}$ and bounded in ${\cal H}$. Then there exists
$\delta_0 >0$ such that, for any $0\leq \delta \leq \delta_0$, for
any $\lambda >0$, the solution (in the sense of Proposition
\ref{propdelta}) ${\bfm u}_\lambda\in \F$ (with $D_t{\bfm u}_\lambda
\in\H)$  of the equation $$\lambda {\bfm u}_\lambda -\delta {\bfm
L}{\bfm u}_\lambda -(1-\delta )\widetilde{{\bfm S}}{\bfm
u}_\lambda-D_t{\bfm u}_\lambda ={\bfm b}_\lambda ,$$ satisfies:
$\exists{\bfm \eta }\in L^2(\Omega )$ such that $\lambda |{\bfm
u}_\lambda|_2^2+|(-\widetilde{{\bfm S}})^{1/2}{\bfm u}_\lambda
-{\bfm \eta}|_2\rightarrow 0$ as $\lambda$ goes to $0$.
\end{proposition}

\noindent {\bf Proof:} Consider the operator $P_\lambda :{\cal
H}\rightarrow {\cal H}$ defined by $P_\lambda ({\bfm b})= ({\bfm
L}-\widetilde{{\bfm S}})(\lambda-\widetilde{{\bfm S}}-D_t
)^{-1}({\bfm b})$. Note that Proposition \ref{propdelta} applies for
all coefficients ${\bfm a}$ and ${\bfm H}$ satisfying Assumption
\ref{hypcontrol}. In particular, it works for ${\bfm
a}=\widetilde{{\bfm a}}$ and ${\bfm H}=0$, so that $P_\lambda$ is
well defined. Lemma \ref{bornex} below proves that
$\|P_\lambda\|_{{\cal H}\to {\cal H}}$ is bounded with a norm that
only depends on the constants $M,C^H_1,C^a_2$ and $C^H_2$ (see
Assumption \ref{hypcontrol}). Therefore, we can choose $\delta_0>0$
such that $\|\delta_0P_\lambda\|_{{\cal H}\to {\cal H}} <1$
(actually $\delta_0<\big[2(2+M+C^H_1)(1+C^a_2+C^H_2)\big]^{-1}$).
For $0<\delta<\delta_0$, we can then define the operator $[{\rm
I}-\delta P_\lambda ]^{-1}:{\cal H}\longrightarrow {\cal H}$. Note
that $(\lambda -\delta {\bfm L}-(1-\delta )\widetilde{{\bfm
S}}-D_t)^{-1}=(\lambda -\widetilde{{\bfm S}}-D_t)^{-1}\big[{\rm
I}-\delta P_\lambda \big]^{-1}$. Thanks to Proposition
\ref{convsdt}, it is sufficient to prove that $\big[{\rm I}-\delta
P_\lambda \big]^{-1}({\bfm b}_\lambda)$ is convergent in $\H_{-1}$.
But $\big[{\rm I}-\delta P_\lambda \big]^{-1}({\bfm
b}_\lambda)=\sum_{n=0}^\infty (\delta P_\lambda)^n({\bfm
b}_\lambda)$. Lemma \ref{bornex} ensures that the sum converges
uniformly with respect to $\lambda>0$. It just remains to prove
that, for each fixed $n\geq 0$, $ ((\delta P_\lambda)^n({\bfm
b}_\lambda))_\lambda$ converges in $\H_{-1}$. This can be proved by
induction on $n\in \nat$. For $n=0$, $({\bfm b}_\lambda
)_{\lambda>0}$ is convergent by assumption. Then, if the family $
((\delta P_\lambda)^n({\bfm b}_\lambda))_\lambda$ is convergent in
$\H_{-1}$, we can apply Proposition \ref{convsdt} to ensure that the
family  $\big((-\widetilde{{\bfm S}})^{1/2}(\lambda-\widetilde{{\bfm
S}}-D_t)^{-1} [(\delta P_\lambda)^n({\bfm b}_\lambda)]\big)_\lambda$
converges in $L^2(\Omega)$. This implies the convergence of
$((\delta P_\lambda)^{n+1}({\bfm b}_\lambda))_\lambda$ in $\H_{-1}$.
\qed

\begin{lemma}\label{bornex}
The norms of $P_\lambda:({\cal H},\|\cdot\|_{-1})\rightarrow
(\H_{-1},\|\cdot\|_{-1})$ and $P_\lambda:({\cal H},\|\cdot\|_{{\cal
H}})\rightarrow ({\cal H},\|\cdot\|_{{\cal H}})$ are both bounded
from above by $2(2+M+C^H_1)(1+C^a_2+C^H_2)$.
\end{lemma}

\vspace{2mm} \noindent {\bf Proof :} Fix ${\bfm b}\in {\cal H}$.
 Let ${\bfm u}_\lambda\in\F$ (with $D_t{\bfm u}_\lambda\in\H$)
be the solution of the equation (apply Proposition \ref{propdelta}
with ${\bfm a}=\widetilde{{\bfm a}}$, ${\bfm H}=0$, ${\bfm h}=0$ and
$m=1$)
\begin{equation*}
 \forall {\bfm \varphi}\in\F,\quad \lambda({\bfm u}_\lambda,{\bfm \varphi})_2+\langle{\bfm u}_\lambda,{\bfm \varphi}\rangle_1-(D_t{\bfm u}_\lambda,
 {\bfm \varphi})_2=\langle{\bfm b},{\bfm \varphi}\rangle_{-1,1}.
\end{equation*}
It derives from \eqref{estimates1} that $\lambda|{\bfm
u}_\lambda|_2^2+\|{\bfm u}_\lambda\|_1^2\leq \|{\bfm b}\|_{-1}^2$,
in such a way that
$$\|P_\lambda({\bfm b})\|_{-1} = \| ({\bfm
L}-\widetilde{{\bfm S}}){\bfm u}_\lambda\|_{-1} \leq
(1+M+C^H_1)\|{\bfm u}_\lambda\|_1\leq (1+M+C^H_1)\|{\bfm
b}\|_{-1}.$$ This proves the first point.\\
Consider now ${\bfm u}\in \F$ with $D_t{\bfm u}\in \H$. An easy
computation proves that, for any $s\in\R^*$ and ${\bfm \varphi}\in
{\cal C}$,\
\begin{equation}\label{calcta}
\begin{split}
{\bfm T}_a\big((-\widetilde{{\bfm S}})^{1/2}{\bfm
u},(-\widetilde{{\bfm S}})^{1/2}\Lambda_s{\bfm \varphi}\big)&
=-\Lambda_{-s}{\bfm T}_a\big((-\widetilde{{\bfm S}})^{1/2}{\bfm
u},(-\widetilde{{\bfm S}})^{1/2}{\bfm \varphi}\big)\\ & -{\bfm
T}_a\big((-\widetilde{{\bfm S}})^{1/2}\Lambda_s{\bfm
u},(-\widetilde{{\bfm S}})^{1/2}T_{s,0}{\bfm \varphi}\big)\\ & \leq
C^a_2\|{\bfm u}\|_1\|{\bfm \varphi}\|_1+M\|D_t{\bfm u}\|_1\|{\bfm
\varphi}\|_1.
\end{split}
\end{equation}
In the above inequalities, we use $\|{\bfm u}\|_1=\|T_{s,0}{\bfm
u}\|_1$ and $\|\Lambda_s{\bfm u}\|_1\leq\|D_t{\bfm u}\|_1$. This
latter point can be proved for ${\bfm u}\in {\cal C}$ as follows
$$\|\Lambda_s{\bfm u}\|_1^2=-(\Lambda_s{\bfm u},\widetilde{{\bfm S}}\Lambda_s{\bfm u})_2=-\int_0^1\int_0^1(
D_tT_{r,0}{\bfm u},\widetilde{{\bfm S}}D_tT_{u,0}{\bfm
u})_2\,dr\,du\leq -(D_t{\bfm u},\widetilde{{\bfm S}}D_t{\bfm u})_2
.$$ The general case is treated by density arguments.\\ As in
\eqref{calcta}, we have ${\bfm T}_H\big((-\widetilde{{\bfm
S}})^{1/2}{\bfm u},(-\widetilde{{\bfm S}})^{1/2}\Lambda_s{\bfm
\varphi}\big)\leq C^H_2\|{\bfm u}\|_1\|{\bfm
\varphi}\|_1+C^H_1\|D_t{\bfm u}\|_1\|{\bfm \varphi}\|_1$. Hence,
$$\|({\bfm L}-\widetilde{{\bfm S}})({\bfm
u})\|_T\leq (C^H_2+C^a_2)\|{\bfm u}\|_1+(C^H_1+M+1)\|D_t{\bfm
u}\|_1.$$ Then, Proposition \ref{propdelta} ensures that $D_t{\bfm
u}_\lambda\in \H$ and  $\|D_t{\bfm u}_\lambda\|_1\leq 2\|{\bfm
b}\|_T+2(C^H_2+C^a_2)\|{\bfm b}\|_{-1}$ (see \eqref{estimates2}) so
that we finally obtain
\begin{equation}\label{inq2}
\begin{split}
\|P_\lambda({\bfm b})\|_{T} & \leq (C^H_2+C^a_2)\|{\bfm
b}\|_{-1}+2(C^H_1+M+1)\big(\|{\bfm b}\|_T+(C^H_2+C^a_2)\|{\bfm
b}\|_{-1}\big).
\end{split}
\end{equation}
The result follows. \qed

\noindent {\bf Proof of Proposition \ref{propcon}:} The last step
before proving Proposition \ref{propcon} consists in lifting the
restriction of the smallness of $\delta_0$. The previous
construction provides us with $\delta_0$ strictly less than $1$. We
perform an induction to get round this restriction whose
initialization is the construction of $\delta_0$. The second step
consists in iterating our arguments to the operator
\begin{multline*}
\lambda -(\delta _0+\delta _1){\bfm L}
-(1-\delta_0-\delta_1)\widetilde{{\bfm S}}-D_t\\ =\big[{\rm
I}-\delta _1({\bfm L}-\widetilde{{\bfm S}})[\lambda -\delta_0 {\bfm
L}-(1-\delta _0)\widetilde{{\bfm S}}-D_t]^{-1}\big](\lambda -\delta
_0{\bfm L}-(1-\delta _0)\widetilde{{\bfm S}}-D_t).
\end{multline*}
We exactly repeat the proof of Proposition \ref{rec} except that the
operator $\lambda -(1-\delta _0-\delta _1)\widetilde{{\bfm
S}}-(\delta_0+\delta _1) {\bfm L}-D_t$ plays the role of the
operator $\lambda -(1-\delta _0)\widetilde{{\bfm S}}-\delta_0 {\bfm
L}-D_t$ and we apply Proposition \ref{rec} with the operator
$\lambda -(1-\delta _1)\widetilde{{\bfm S}}-\delta_1 {\bfm L}-D_t$
instead of applying Proposition \ref{convsdt} with $ \lambda
-\widetilde{{\bfm S}}-D_t$. Of course, a restriction about the
smallness of $\delta_1$ is imposed by this procedure. Even if it
means substituting $\widetilde{a}$ with $m\widetilde{a}$, we assume,
without loss of generality, that $m=1$. Thus  Lemma \ref{bornex}
remains valid for the operator $P^1_\lambda:{\cal H}\rightarrow
{\cal H} $ defined by $P^1_\lambda({\bfm b})=({\bfm
L}-\widetilde{{\bfm S}})(\lambda-(1-\delta_0)\widetilde{{\bfm
S}}-\delta_0{\bfm L}-D_t )^{-1}({\bfm b})$. This is of the utmost
importance because that means that we can choose
$\delta_1=\delta_0$. Thus we can iterate these arguments until we
find $\delta _N$ such that $\delta _0+\delta _1+\dots+\delta _N>1$
and such that Proposition \ref{rec} still holds except that $\delta
<\delta_0$ is everywhere replaced by $\delta<\delta_0+\delta
_1+\dots+\delta _N$. Proposition \ref{propcon} follows.\qed

Now let us prove that the drift ${\bfm b}$ of the diffusion process
$X$ fulfills the assumptions of Proposition \ref{propcon}. To this
purpose, let us establish

\begin{lemma}\label{lemmeb}
For each $i\in\{1,\dots,d\}$, ${\bfm b}_i$ belongs to $\H_{-1}$ and
$\forall s\in\R,\forall {\bfm \varphi }\in {\cal C}$, $$\langle
{\bfm b}_i,\Lambda_s{\bfm \varphi }\rangle_{-1,1} \leq
(C^a_2+C^H_2)|(\widetilde{\bfm{a}}E_i,E_i)_2|^{1/2}\|{\bfm \varphi
}\|_1.$$
\end{lemma}

\noindent {\bf Proof:} Let $(E_1,\dots,E_d)$ be the canonical basis
of $\reel ^d$. Then we have
\begin{eqnarray*}
({\bfm b}_i,{\bfm \varphi })_2 & = & 1/2\sum_{j} \big(e^{2 {\bfm
V}}D_j(e^{-2 {\bfm V}}[\bfm{a}+\bfm{H}]_{ij}), {\bfm \varphi
}\big)_2\\ & =& -1/2\big([\bfm{a}-\bfm{H}] D\bfm{\varphi },
E_i\big)_2
\\& \leq & 1/2\big|\big(\bfm{a} D\bfm{\varphi },
E_i\big)_2\big|+1/2\big|\big(\bfm{H} D\bfm{\varphi },
E_i\big)_2 \big|\\
 & \stackrel{Cauchy-Schwarz}{\leq} & M\|\bfm{\varphi
 }\|_1|(\widetilde{\bfm{a}}E_i,E_i)_2|^{1/2} +C^H_1\|\bfm{\varphi
 }\|_1|(\widetilde{\bfm{a}}E_i,E_i)_2|^{1/2}
\end{eqnarray*}
and this proves the first point. Then, $\forall s>0,\forall {\bfm
\varphi }\in {\cal C}$, we have
\begin{eqnarray*}
\langle {\bfm b}_i,\Lambda_s{\bfm \varphi }\rangle_{-1,1} & =&
-(1/2)\big([{\bfm a}+{\bfm H}]E_i,\Lambda_sD{\bfm \varphi }\big)_2\\
 & = & (1/2)\big(\Lambda_{-s}[{\bfm a}+{\bfm H}]E_i,D{\bfm \varphi }\big)_2\\ & \stackrel{\text{Assumption }\ref{hypcontrol} }{\leq}
  & (C^a_2+C^H_2)|(\widetilde{\bfm{a}}E_i,E_i)_2|^{1/2}\|{\bfm \varphi }\|_1\qed
\end{eqnarray*}

%%%%%%%%%%%%%%%%%%%%%%%%%%%%%%%%%%%%%%%%%%%%%%%%%%%%%%%%%%%%%%%%%%%%%%%%%%%%%%%%%%%%%%%%%%%%%%%%%%%%%%%%%%%%%%

\section{Itô's formula}\label{ito}

%%%%%%%%%%%%%%%%%%%%%%%%%%%%%%%%%%%%%%%%%%%%%%%%%%%%%%%%%%%%%%%%%%%%%%%%%%%%%%%%%%%%%%%%%%%%%%%%%%%%%%%%%%%%%
We are not in a lucky situation of working on an explicit Dirichlet
form connected with the generator in $L^2(\Omega,\pi)$ of $Y$,
wrongly denoted by $[{\bfm L}+D_t]$. This raises the following
issue: given a function ${\bfm f}\in L^2(\Omega)$ and the function
${\bfm u}_\lambda$ that weakly solves (see Proposition
\ref{propdelta}) $\lambda {\bfm u}_\lambda -({\bfm L}+D_t){\bfm
u}_\lambda={\bfm f}$, does the "Ito formula" apply to ${\bfm
u}_\lambda$ and to the process $Y$. Indeed, it is not clear that the
construction of ${\bfm u}_\lambda$ in Proposition \ref{propdelta}
belongs to the domain of the generator of $Y$. The key tool is the
regular approximation $({\bfm u}_{\lambda,\delta})_\delta$ provided
by Proposition
\ref{propdelta} for a suitable function ${\bfm f}$.\\
Let us consider a standard $1$-dimensional Brownian motion
$\{B'_t;t\geq 0\}$ that is independent of $\{B_t;t\geq 0\}$ in such
a way that $\{(B'_t,B_t);t\geq 0\}$ is a standard $d+1$-dimensional
Brownian motion. Define then the $d+1$-dimensional diffusion process
$X^{\omega,\delta}$, starting from $0$, as the solution of the SDE:
\begin{equation}\label{eq_diffdelta}
 X^{\omega,\delta}_t=\int_0^t\left[ \begin{array}{c}
                                 1 \\
                                 b(X^{\omega,\delta}_r,\omega) \\
                               \end{array}\right]\,dr+\int_0^t\left[ \begin{array}{cc}
                                 \sqrt{\delta}  & 0 \\
                                 0 & \sigma(X^{\omega,\delta}_r,\omega) \\
                               \end{array}\right]\,d(B',B)_r.
\end{equation}
The associated diffusion in random medium  $Y^\delta$ defined by
$Y^\delta_t(\omega)=\tau_{X^{\omega,\delta}_t}\omega$ is a
$\Omega$-valued Markov process, which admits $\pi$ as invariant
measure (similar to section \ref{sec_env}). It also defines a
continuous semi-group on $L^2(\Omega)$. The associated
(non-symmetric) Dirichlet form is given by
$\eqref{def_dirichletform}$ (with $\theta=1$) with domain
$\F\times\F$ and satisfies a weak sector condition (see \cite[Ch. 1,
Sect 2.]{ma} for the definition). The generator ${\bfm L}^\delta$ is
defined on ${\rm Dom}({\bfm L}^\delta)=\{{\bfm u}\in\F;
B_{\lambda,\delta}({\bfm u},\cdot)\text{ is
}L^2(\Omega)\text{-continuous}\}$ (see \cite[Ch. 1, Sect 2.]{ma} for
further details). It coincides on ${\cal C}$ with $ {\bfm
L}+D_t+(\delta/2)D^2_t$. Since $b$ and $\sigma$ are globally
Lipschitz (Assumption \ref{hypregularity}), classical tools of SDE
theory ensures that
\begin{equation}\label{convxdelta}
 \int_\Omega\E\big[\sup_{0\leq t \leq T}|(t,X^\omega_t)-X_t^{\omega,\delta}|^2\big]\,d\pi\rightarrow 0 \text{ as }\delta\text{ goes to }0,
\end{equation}
where both diffusions start from $0$.

\begin{proposition}\label{prop_ito}
Let ${\bfm f}\in L^2(\Omega)$ and a family $({\bfm
u}_\lambda)_{\lambda>0}$ in $\F$ such that:

\noindent 1) $\forall {\bfm \varphi}\in \F$, $B_\lambda({\bfm
u}_\lambda,{\bfm \varphi})=({\bfm f},{\bfm \varphi})_2 $,

\noindent 2) for each $\lambda>0$, there exists a sequence $({\bfm
u}_{\lambda,\delta})_{\delta>0}$ in  $\F$ that converges in $\H$
towards ${\bfm u}_\lambda$. Moreover $({\bfm
u}_{\lambda,\delta})_{\delta>0}\in{\rm Dom}({\bfm L}^\delta)$ and
satisfies $\lambda {\bfm u}_{\lambda,\delta} -{\bfm L}^\delta{\bfm
u}_{\lambda,\delta}={\bfm f}$.

\noindent 3) for each fixed $\lambda>0$, $(D_t{\bfm
u}_{\lambda,\delta})_\delta$ is bounded in $L^2(\Omega)$.

\noindent 4) each function ${\bfm u}_{\lambda,\delta}$ has
continuous trajectories, that is, for $\mu$ almost every $\omega\in
\Omega$, the function $(t,x)\in\R^{d+1}\mapsto {\bfm
u}_{\lambda,\delta}(\tau_{t,x}\omega)$ is continuous.\\
Then, $\P_\pi$ a.s., the following formula holds
\begin{equation*}
{\bfm u}_\lambda(Y_t)={\bfm
u}_\lambda(Y_0)+\int_{0}^{t}(\lambda{\bfm u}_\lambda-{\bfm f})
(Y_r)\,dr+\int_{0}^{t}\nabla ^\sigma{\bfm u}_\lambda^*(Y_r)\,dB_r
\end{equation*}
where $\prob_\pi $ is the law of the process $Y$ starting with
initial distribution $\pi $ on $\Omega $.
\end{proposition}

\noindent {\bf Proof:} Since ${\bfm u}_{\lambda,\delta}\in {\rm
Dom}({\bfm L}^\delta)$ and $\lambda{\bfm u}_{\lambda,\delta}-{\bfm
L}^\delta{\bfm u}_{\lambda,\delta}={\bfm f}$, we can write (see
Lemma \ref{itoreg} below)
\begin{equation}\label{itoudelta}
\begin{split}
{\bfm u}_{\lambda,\delta}(Y^\delta_t) & -{\bfm
u}_{\lambda,\delta}(Y^\delta_0) \\ & =\int_0^t{\bfm L}^\delta{\bfm
u}_{\lambda,\delta}(Y^\delta_r)\,dr+\delta^{1/2}\int_0^tD_t{\bfm
u}_{\lambda,\delta}(Y^\delta_r)\,dB'_r+\int_0^t\nabla^\sigma{\bfm
u}_{\lambda,\delta}^*(Y^\delta_r)\,dB_r\\ &  = \int_0^t[\lambda{\bfm
u}_{\lambda,\delta}-{\bfm
f}](Y^\delta_r)\,dr+\delta^{1/2}\int_0^tD_t{\bfm
u}_{\lambda,\delta}(Y^\delta_r)\,dB'_r+\int_0^t\nabla^\sigma{\bfm
u}_{\lambda,\delta}^*(Y^\delta_r)\,dB_r.
\end{split}
\end{equation}
Thanks to \eqref{convxdelta}, the convergence, as $\delta\to 0$, of
$({\bfm u}_{\lambda,\delta})_{\lambda,\delta}$ towards ${\bfm
u}_{\lambda}$ in $\H$ and the boundedness of $(D_t{\bfm
u}_{\lambda,\delta})_\delta$ in $L^2(\Omega)$, we can pass to the
limit in \eqref{itoudelta} and complete the proof.
 \qed

\begin{lemma}\label{itoreg}
Keeping the notations of Proposition \ref{prop_ito}, the following
formula holds, $\P_\pi$ a.s.,
$${\bfm u}_{\lambda,\delta}(Y^\delta_t)-{\bfm
u}_{\lambda,\delta}(Y^\delta_0)  =\int_0^t{\bfm L}^\delta{\bfm
u}_{\lambda,\delta}(Y^\delta_r)\,dr+\delta^{1/2}\int_0^tD_t{\bfm
u}_{\lambda,\delta}(Y^\delta_r)\,dB'_r+\int_0^t\nabla^\sigma{\bfm
u}_{\lambda,\delta}^*(Y^\delta_r)\,dB_r.$$
\end{lemma}
\noindent {\bf Proof:} Since ${\bfm u}_{\lambda,\delta}\in {\rm
Dom}({\bfm L}^\delta)$ , the difference ${\bfm
u}_{\lambda,\delta}(Y^\delta_t)-{\bfm
u}_{\lambda,\delta}(Y^\delta_0)-\int_0^t{\bfm L}^\delta{\bfm
u}_{\lambda,\delta}(Y^\delta_r)\,dr$ is a square-integrable
continuous $\P_\pi$-martingale, denoted by $M^\delta_t$. Moreover,
for a function ${\bfm \varphi}\in {\cal C}$, the classical Ito
formula yields ${\bfm \varphi}(Y^\delta_t)-{\bfm
\varphi}(Y^\delta_0) =\int_0^t{\bfm L}^\delta{\bfm
\varphi}(Y^\delta_r)\,dr+\delta^{1/2}\int_0^tD_t{\bfm
\varphi}(Y^\delta_r)\,dB'_r+\int_0^t\nabla^\sigma{\bfm
\varphi}^*(Y^\delta_r)\,dB_r $. Then the process $t\mapsto {\bfm
u}_{\lambda,\delta}(Y^\delta_t)-{\bfm \varphi}(Y^\delta_t)$ is a
continuous semimartingale and Theorem 32 in \cite[Ch. 2, Sect.
7]{protter} (applied with the function $x\in \R\mapsto x^2$) yields
$\P_\pi$ a.s.,
\begin{equation}
\begin{split}
({\bfm u}_{\lambda,\delta}(Y^\delta_t)- & {\bfm
\varphi}(Y^\delta_t))^2\\ & =({\bfm
u}_{\lambda,\delta}(Y^\delta_t)-{\bfm
\varphi}(Y^\delta_0))^2+2\int_0^t({\bfm u}_{\lambda,\delta}-{\bfm
\varphi}){\bfm L}^\delta({\bfm u}_{\lambda,\delta}-{\bfm
\varphi})(Y^\delta_r)\,dr\\ & +2\int_0^t({\bfm
u}_{\lambda,\delta}-{\bfm
\varphi})(Y^\delta_r)\,\big(dM^\delta_r-\delta^{1/2}D_t{\bfm
\varphi}(Y^\delta_r)\,dB'_r-\nabla^\sigma{\bfm
\varphi}^*(Y^\delta_r)\,dB_r\big)\\ &
+2\big[M-\int_0^\cdot\delta^{1/2}D_t{\bfm
\varphi}(Y^\delta_r)\,dB'_r-\int_0^\cdot\nabla^\sigma{\bfm
\varphi}^*(Y^\delta_r)\,dB_r\big]_t,
\end{split}
\end{equation}
where $[X]$ stands for the quadratic variations of the martingale
$X$. Integrating with respect to the measure $\pi$, the martingale
term vanishes and we deduce
\begin{equation}\label{estapprox}
 \E_\pi\big(2\big[M-\int_0^\cdot\delta^{1/2}D_t{\bfm
\varphi}(Y^\delta_r)\,dB'_r-\int_0^\cdot\nabla^\sigma{\bfm
\varphi}^*(Y^\delta_r)\,dB_r\big]_t\big)\leq
2B_{\lambda,\delta}({\bfm u}_{\lambda,\delta}-{\bfm \varphi},{\bfm
u}_{\lambda,\delta}-{\bfm \varphi}).
\end{equation}
Choosing a sequence $({\bfm \varphi}_n)_n$ in ${\cal C}$ that
converges in $\F$ towards $ {\bfm u}_{\lambda,\delta}$, we easily
complete the proof with the help of \eqref{estapprox}.\qed

 Note that
the time reversed process $t \mapsto Y_{T-t}^\delta$ is a Markov
process with respect to the backward filtration $({\cal G}^\delta
_t)_{0\leq t\leq T}$, where ${\cal G}^\delta_s$ is the
$\sigma$-algebra on $\Omega $ generated by $\left\{ Y^\delta_r;t\leq
r\leq T\right\}$, and admits the adjoint operator $({\bfm L}^\delta
)^*$ of ${\bfm L}^\delta$ in $L^2(\Omega,\pi)$ as generator, which
coincides on ${\cal C}$ with $ {\bfm L}^*-D_t+(\delta/2)D^2_t$. From
\eqref{convxdelta}, $t \mapsto Y_{T-t}^\delta$ approximates the
process $ t \mapsto Y_{T-t}$ as $\delta$ tends to $0$. It is then
readily seen that we can slightly modify the proof of Proposition
\ref{prop_ito} and prove the
\begin{proposition}\label{prop_ito2}
Let ${\bfm f}\in L^2(\Omega)$ and a family $({\bfm
u}_\lambda)_{\lambda>0}$ in $\F$ such that:

\noindent 1) $\forall {\bfm \varphi}\in \F$, $B_\lambda({\bfm
\varphi},{\bfm u}_\lambda)=({\bfm f},{\bfm \varphi})_2 $,

\noindent 2) for each $\lambda>0$, there exists a sequence $({\bfm
u}_{\lambda,\delta})_{\delta>0}$ in  $\F$ that converges in $\H$
towards ${\bfm u}_\lambda$. Moreover $({\bfm
u}_{\lambda,\delta})_{\delta>0}\in{\rm Dom}({\bfm L}^\delta)^*$ and
satisfies $\lambda {\bfm u}_{\lambda,\delta} -({\bfm
L}^\delta)^*{\bfm u}_{\lambda,\delta}={\bfm f}$.

\noindent 3) for each fixed $\lambda>0$, $(D_t{\bfm
u}_{\lambda,\delta})_\delta$ is bounded in $L^2(\Omega)$.

\noindent 4) each function ${\bfm u}_{\lambda,\delta}$ has
continuous trajectories, that is, for $\mu$ almost every $\omega\in
\Omega$, the function $(t,x)\in\R^{d+1}\mapsto {\bfm
u}_{\lambda,\delta}(\tau_{t,x}\omega)$ is continuous.\\
Then, $\P_\pi$ a.s., the following formula holds
\begin{equation*}
{\bfm u}_\lambda(Y_T-t)={\bfm
u}_\lambda(Y_T)+\int_{0}^{t}(\lambda{\bfm u}_\lambda-{\bfm f})
(Y_{T-r})\,dr+(M_t-M_0)
\end{equation*}
where $M$ is a martingale with respect to the backward filtration
$({\cal G} _t)_{0\leq t\leq T}$, and ${\cal G}_s$ is the
$\sigma$-algebra on $\Omega $ generated by $\left\{ Y_r;t\leq r\leq
T\right\}$. Moreover, the quadratic variations of $M$ exactly match
$\int_{0}^{t}\nabla ^\sigma{\bfm u}_\lambda^*\cdot\nabla
^\sigma{\bfm u}_\lambda(Y_{T-r})\,dr$.
\end{proposition}

%%%%%%%%%%%%%%%%%%%%%%%%%%%%%%%%%%%%%%%%%%%%%%%%%%%%%%%%%%%%%%%%%%%%%%%%%%%%%%%%%%%%%%%%%%

\section{Ergodic Theorem}\label{sec_ergodic}

%%%%%%%%%%%%%%%%%%%%%%%%%%%%%%%%%%%%%%%%%%%%%%%%%%%%%%%%%%%%%%%%%%%%%%%%%%%%%%%%%%%%%%%%%%
Let us now exploit the ergodic properties of the operator
$\widetilde{{\bfm S}}$ stated in Assumption \ref{ergodicity} and
prove
\begin{theorem}\label{theo_erg}
Let ${\bfm f}\in L^1(\Omega)$. Then
$$\E_\pi\Big|\frac{1}{t}\int_0^t{\bfm f}(Y_r)\,dr-\pi({\bfm f})\Big|\rightarrow 0 \text{ as } t \text{ goes to } \infty.$$
\end{theorem}

\noindent {\bf Proof:} We suppose at first that ${\bfm f}\in {\cal
C}$. Even if it means considering ${\bfm f}-\pi({\bfm f})$ instead
of ${\bfm f}$, we assume that $\pi({\bfm f})=0$. Clearly, ${\bfm
f}\in {\rm Dom}(D_t)$ and Proposition \ref{propdelta} applies. For
each $\lambda>0$, it provides us with a function ${\bfm
u}_\lambda\in\F$ such that
\begin{equation}\label{eq_erg}
\forall{\bfm \varphi}\in \F,\quad B_\lambda({\bfm u}_\lambda,{\bfm
\varphi})=({\bfm f},{\bfm \varphi})_2.
\end{equation}
Moreover, \eqref{estimates1} and \eqref{estimates2} ensures that the
families $(\lambda{\bfm u}_\lambda)_\lambda$, $(\lambda D_t{\bfm
u}_\lambda)_\lambda $ and $(\lambda^{1/2}(-\widetilde{{\bfm
S}})^{1/2}{\bfm u}_\lambda)_\lambda$ are bounded in $L^2(\Omega)$.
Even if it means considering a subsequence, we assume that
$(\lambda{\bfm u}_\lambda)_\lambda$, $(\lambda D_t{\bfm
u}_\lambda)_\lambda$ and $(\lambda^{1/2}(-\widetilde{{\bfm
S}})^{1/2}{\bfm u}_\lambda)_\lambda$ weakly converge respectively to
${\bfm g}$, ${\bfm g}'$ and ${\bfm G}$ in $L^2(\Omega)$. Since the
operator $D_t$ is closed, it turns out that ${\bfm g}'=D_t{\bfm g}
$. Let us now prove now that ${\bfm g}\in{\rm Dom}({\bfm L})$.
Consider ${\bfm \varphi}\in{\rm Dom }({\bfm L}^*)$ . Then we derive
from \eqref{eq_erg} that
\begin{equation*}
\begin{split}
\lambda({\bfm f},{\bfm \varphi})_2=\lambda B_\lambda({\bfm
u}_\lambda,{\bfm \varphi})=\lambda^2( {\bfm u}_\lambda,{\bfm
\varphi})_2 -(\lambda {\bfm u}_\lambda,{\bfm L}^*{\bfm
\varphi})_2-(\lambda D_t{\bfm u}_\lambda,{\bfm \varphi})_2.
\end{split}
\end{equation*}
Passing to the limit as $\lambda$ goes to $0$, we deduce $({\bfm
g},{\bfm L}^*{\bfm \varphi})_2=-(D_t{\bfm g},{\bfm \varphi})_2$.
Hence ${\bfm g}\in {\rm Dom}({\bfm L}^{**})={\rm Dom}({\bfm
L})\subset \H$ and ${\bfm L}{\bfm g}=-D_t{\bfm g}$. In particular
$$m\|{\bfm g}\|_1^2\leq -({\bfm g},{\bfm L}{\bfm g})_2=(D_t{\bfm
g},{\bfm g})_2=0$$ so that ${\bfm g}\in{\rm Dom}((-\widetilde{{\bfm
S}})^{1/2})$ and $(-\widetilde{{\bfm S}})^{1/2}{\bfm g}=0$. As a
consequence, ${\bfm g}\in {\rm Dom}( \widetilde{{\bfm S}})$ and
$\widetilde{{\bfm S}}{\bfm g}=0$. From Assumption \ref{ergodicity},
${\bfm g}$ is invariant under space translations in such a way that
$D_t{\bfm g}=-{\bfm L}{\bfm g}=0$ and ${\bfm g}$ is also invariant
under time translations. Thus the ergodicity of the measure $\mu$
implies that ${\bfm g}$ is constant ($\mu$ a.s.). Choosing ${\bfm
\varphi}$ equal to the constant function $\one$ in \eqref{eq_erg},
we deduce ${\bfm g}=0$. We now aim at proving that the convergence
of $(\lambda{\bfm u}_\lambda)_\lambda $ towards $0$ holds in the
strong sense. In what follows, we make no distinction between
$0\in\R$ and the constant function that matches $0$ over $\Omega$.
We just have to write
\begin{equation*}
0=(0,{\bfm f})_2=\lim_{\lambda\to 0}(\lambda{\bfm u}_\lambda,{\bfm
f})_2 = \lim_{\lambda\to 0}B_\lambda(\lambda{\bfm
u}_\lambda,\lambda{\bfm u}_\lambda)_2\geq\limsup_{\lambda\to
0}|\lambda{\bfm u}_\lambda|_2^2.
\end{equation*}
Note now that the approximating family $({\bfm
u}_{\lambda,\delta})_\delta $ provided by Proposition
\ref{propdelta} is given by ${\bfm
u}_{\lambda,\delta}(\omega)=\int_0^\infty e^{-\lambda
r}\E_0[f(X^{\omega,\delta}_r,\omega)]\,dr $. For each $(t,x)\in
\R^{d+1}$, the law of the process
$(t,x)+X^{\tau_{t,x}\omega,\delta}$, $X^{\tau_{t,x}\omega,\delta}$
starting from $0\in\R^{d+1}$, is the same as the law of the process
$X^{\omega,\delta}$ starting from $(t,x)\in\R^{d+1}$ (see the proof
at the end of Section \ref{sec_inv}). Hence ${\bfm
u}_{\lambda,\delta}(\tau_{t,x}\omega)=\int_0^\infty e^{-\lambda
r}\E_{t,x}[f(X^{\omega,\delta}_r,\omega)]\,dr $. Since ${\bfm f}$ is
smooth and $X^{\omega,\delta}$ is a Feller process, ${\bfm
u}_{\lambda,\delta}$ has continuous trajectories. Thus Proposition
\ref{prop_ito} applies and it yields
\begin{equation*}
 \int_0^t{\bfm
f}(Y_r)\,dr =({\bfm u}_{\lambda}(Y^\delta_0)-{\bfm
u}_{\lambda}(Y_t))+\int_0^t\lambda{\bfm
u}_{\lambda}(Y_r)\,dr+\int_0^t\nabla^\sigma{\bfm
u}_{\lambda}^*(Y_r)\,dB_r.
\end{equation*}
Thanks to \eqref{estimates1} and the invariance of the measure $\pi$
for the process $Y$, we can find a constant $C$, which depends
neither on $\lambda$ nor on $t$, such that
\begin{equation*}
\E_\pi\big| \frac{1}{t}\int_0^t{\bfm f}(Y_r)\,dr\big|^2 \leq
C/(t\lambda)^2+C|\lambda{\bfm u}_{\lambda}|^2_2+C/(t\lambda^{1/2}).
\end{equation*}
It just remains to choose $\lambda$ small enough and then $t$ large
enough to complete the proof in the case ${\bfm f}\in {\cal C}$. The
general case is treated with the density of ${\cal C}$ in
$L^1(\Omega)$ and the invariance of the measure $\pi$. Since it
raises no particular difficulty, details are left to the reader.
\qed

%%%%%%%%%%%%%%%%%%%%%%%%%%%%%%%%%%%%%%%%%%%%%%%%%%%%%%%%%%%%%%%%%%%%%%%%%%%%%%%%%%%%%%%%%%

\section{Invariance principle}\label{sec_inv}

%%%%%%%%%%%%%%%%%%%%%%%%%%%%%%%%%%%%%%%%%%%%%%%%%%%%%%%%%%%%%%%%%%%%%%%%%%%%%%%%%%%%%%%%%%
\noindent {\bf Notation :}\\
\noindent Up to the end of this paper, for $i\in\{1,\dots,d\}$ we
denote by ${\bfm u}^i_\lambda$ the solution of the equation (in the
weak sense of Proposition \ref{propdelta})
$$\lambda  {\bfm u}^i_\lambda-{\bfm L}{\bfm u}^i_\lambda-D_t{\bfm
u}^i_\lambda={\bfm b}_i.$$ From Proposition \ref{propcon}, there
exists ${\bfm \xi }_i\in (L(\Omega))^d$ such that $\lambda |{\bfm
u}^i_\lambda|_2^2+|{\bfm \xi }_i-\nabla ^\sigma {\bfm
u}^i_\lambda|_2\rightarrow 0$ as $\lambda$ goes to $0$.\qed
\\Applying the Ito formula (see Proposition \ref{prop_ito}) to the
function ${\bfm u}_{\varepsilon ^2}$ yields
$$\varepsilon X^{\omega }_{t/\varepsilon^2}=H^{\varepsilon,\omega } _t+\varepsilon
\int_{0}^{t/\varepsilon^2}(\sigma  +\nabla ^\sigma
u_\lambda^*)(r,X^{\omega }_r,\omega)\,dB_r,$$ where
$$H^{\varepsilon,\omega } _t=\varepsilon ^3\int_{0}^{t/\varepsilon
^2}u_{\varepsilon ^2}(r,X^{\omega }_r,\omega )\,dr-\varepsilon
u_{\varepsilon ^2}(t/\varepsilon ^2,X^{\omega }_{t/\varepsilon
^2},\omega )+\varepsilon u_{\varepsilon ^2}(0,0,\omega).$$ For the
reader's convenience, it is worth recalling that $Y_t=\tau
_{t,X^{\omega }_t}$ and $\prob_\pi $ is the law of the process $Y$
with initial distribution $\pi $. We want to show that the finite
dimensional distributions of the process $H^{\varepsilon,\omega }$
converges in $\prob_\pi $-probability to $0$. Using the
Cauchy-Scharz inequality and the invariance of the measure $\pi $,
we get the estimate
$$\E_{\pi} [(H^{\varepsilon,\omega } _t)^2] \leq
3(2+t^2)\varepsilon ^2|u_{\varepsilon^2}|^2_{2}$$ and this latter
quantity converges to
$0$ as $\varepsilon $ goes to $0$.\\
Let us now investigate the convergence of the process $t\mapsto
\varepsilon \int_{0}^{t/\varepsilon^2}(\sigma  +\nabla ^\sigma
u_\lambda^*)(Y_r)\,dB_r$ whose quadratic variations are given by
\begin{equation*}
\begin{split}
\varepsilon ^2\int_{0}^{t/\varepsilon ^2}&({\bfm \sigma} +\nabla
^\sigma {\bfm u}_{\varepsilon ^2}^*) ({\bfm \sigma }+\nabla ^\sigma
{\bfm u}_{\varepsilon ^2}^*)^*(Y_r)\,dr = \varepsilon
^2\int_{0}^{t/\varepsilon ^2}( {\bfm \sigma} +{\bfm \xi }^*)({\bfm
\sigma }+{\bfm \xi }^*)^*(Y_r)\,dr\\&+\Big( \varepsilon
^2\int_{0}^{t/\varepsilon ^2}({\bfm \sigma} +\nabla ^\sigma {\bfm
u}_{\varepsilon ^2}^*) ({\bfm \sigma }+\nabla ^\sigma {\bfm
u}_{\varepsilon ^2}^*)^*(Y_r)\,dr - \varepsilon
^2\int_{0}^{t/\varepsilon ^2}({\bfm \sigma }+{\bfm \xi }^*)({\bfm
\sigma} +{\bfm \xi }^*)^*(Y_r)\,dr\Big).
\end{split}
\end{equation*}
With the help of Theorem \ref{theo_erg}, the finite dimensional
distributions of the  former term in the right-hand side converge in
$L^1(\P_\pi)$ to the ones of the process $t\mapsto At$, where the
matrix $A$ is given  by
\begin{equation}\label{def_A}
  A=\int_{\Omega }({\bfm \sigma}+{\bfm \xi }^*)({\bfm \sigma}+{\bfm \xi }^*)^*\,d\pi .
\end{equation}
The finite dimensional distributions of the latter term in the
right-hand side converge in $L^1(\P_\pi)$ to $0$. Indeed, after
integrating with respect to the probability measure $\prob_\pi$, it
is bounded by $Ct|\nabla ^\sigma~{\bfm u}_{\varepsilon ^2}-~{\bfm
\xi }|_2^2$. Hence we conclude by applying the central limit theorem
for martingales that the finite dimensional distributions of the
process $\varepsilon X^{\omega} _{t/\varepsilon ^2}$ converge in law
to the ones of the process $A^{1/2}B_t$.

\begin{proposition}
The process $\varepsilon X^{\omega} _{t/\varepsilon ^2}$ is tight in the space $C([0,T];\reel^d)$. Hence it converges in law in the space $C([0,T];\reel^d)$ towards the process $A^{1/2}B_t$.
\end{proposition}

\vspace{2mm} \noindent {\bf Proof :} The next section is devoted to
the proof of the tightness \qed \\Let us now to determine the limit
when the starting point is not $0$ but $x/\varepsilon $.
\begin{equation*}
\begin{split}
\E_{x/\varepsilon }\big[f(\varepsilon X^{\omega }_{t/\varepsilon
^2})\big] & = \E_0\big[f(x +\varepsilon X^{ \tau_{(0,x/\varepsilon
)}\omega }_{t/\varepsilon ^2})\big] \stackrel{\text{in law w.r.t.
}\mu }{=}  \E_0\big[f(x+\varepsilon X^{\omega }_{t/\varepsilon
^2})\big]\\&  \xrightarrow[\varepsilon \rightarrow 0]{\pi \text{
prob}} \E\big[f(x+A^{1/2}B_t)\big]
\end{split}
\end{equation*}
For the first above equality we used the following fact. If
$$X_t=x+\int_{0}^{t}b\left(r, X_r,\omega \right)
\,dr+\int_{0}^{t}\sigma\left(r, X_r,\omega \right)\,dB_r$$
and $Z_t\stackrel{ \Delta }{=}X_t-x$
then $Z_t$ solves the SDE
$$Z_t=\int_{0}^{t}b\left(r, Z_r,\tau_{(0,x)}
\omega \right)\,dr+\int_{0}^{t}\sigma\left(r,Z_r,\tau_{(0,x )}
\omega \right)\,dB_r,$$ so that the law of the process $X^\omega$
starting from $x\in\R^d$ is equal to the law of the process
$x+X^{\tau_x\omega} $ where $X^{\tau_x\omega}$ is starting from $0$.
We sum up:
\begin{theorem}
Let $ f$ be a continuous, bounded function on $\R^d$. Then the
solution $z(x ,t,\omega )$ of the partial differential equation
\eqref{equationeps} with initial condition $z(0,x,\omega)=f(x)$
satisfies the following convergence: $z(x/\varepsilon ,t/\varepsilon
^2,\omega )$ converges in $\pi $-probability as $\varepsilon
\rightarrow 0$ to $\E\left[f(x+A^{1/2}B_t)\right]$, which is the
viscosity solution of the deterministic equation \eqref{eq_lim} with
the same initial condition. The matrix $A$ is given by
$$A=\int_{\Omega }({\bfm \sigma}+{\bfm \xi }^*)({\bfm \sigma}+{\bfm
\xi }^*)^*\,d\pi .$$
\end{theorem}

%%%%%%%%%%%%%%%%%%%%%%%%%%%%%%%%%%%%%%%%%%%%%%%%%%%%%%%%%%%%%%%%%%%%%%%%%%%%%%%%%%%%%%%%%%%%%%%%%%%%%%%%%%%%%%%%%%
\section{Tightness}
%%%%%%%%%%%%%%%%%%%%%%%%%%%%%%%%%%%%%%%%%%%%%%%%%%%%%%%%%%%%%%%%%%%%%%%%%%%%%%%%%%%%%%%%%%%%%%%%%%%%%%%%%%%%%%%%%%

Let us now investigate the tightness in $C([0,T];\reel^d)$ of the
process
$$\varepsilon X^{\omega} _{t/\varepsilon ^2}=\varepsilon
\int_0^{t/\varepsilon ^2}b(r,X^{\omega }_r,\omega )\,dr+\varepsilon
\int_0^{t/\varepsilon ^2}\sigma (r,X^{\omega }_r,\omega )\,dB_r.$$
The tightness of the first term in the above right-hand side is
readily derived from the Burkholder-Davis-Gundy inequality and the
boundedness of the diffusion coefficient ${\bfm \sigma} $.
Concerning the second term, we are going to exploit ideas of
\cite{sethu}  or \cite{wu}.

For any $i\in\{1,\dots,d\}$ and $\lambda >0$, we put ${\bfm
w}_\lambda =(\lambda -{\bfm S})^{-1}{\bfm b}_i\in \H\cap {\rm
Dom}({\bfm S})$
 (see Proposition \ref{existnotime}). Proposition \ref{propdelta} (with $\theta=0$ and ${\bfm H}=0$) also ensures
 that ${\bfm
w}_\lambda\in\F$, $D_t{\bfm w}_\lambda\in\H$. For each fixed
$\lambda>0$,
  we can find a sequence $({\bfm \psi}_\lambda^n)_n$ in ${\cal
  C}$ such that $\|{\bfm \psi}_\lambda^n-{\bfm w}_\lambda\|_1+\|D_t{\bfm \psi}_\lambda^n- D_t{\bfm
w}_\lambda\|_1$ converges to $0$ as $n$ goes
  to $\infty$. Define ${\bfm A}^n_\lambda= (1/2)\sum_{k,l}D_l\big({\bfm H}_{kl}D_k{\bfm
  \psi}_\lambda^n\big)$. From Proposition \ref{propdelta}, we can find two
  sequences $(\overline{{\bfm v}}_\lambda^n)_n\subset \F\cap {\rm Dom}({\bfm L})$ and $(\underline{{\bfm v}}_\lambda^n)_n\subset \F\cap {\rm Dom}({\bfm
  L}^*)$ that respectively solve the equations $(\lambda-{\bfm L})\overline{{\bfm v}}_\lambda^n={\bfm
  b}_i-{\bfm A}^n_\lambda$ and $(\lambda-{\bfm L}^*)\underline{{\bfm v}}_\lambda^n={\bfm
  b}_i+{\bfm A}^n_\lambda$. Moreover, the functions $\overline{{\bfm
  v}}_\lambda^n$ and $\underline{{\bfm v}}_\lambda^n$ possess a
  corresponding approximation sequence $(\overline{{\bfm v}}_{\lambda,\delta}^n)_{\delta>0} $ and $(\underline{{\bfm v}}_{\lambda,\delta}^n)_{\delta>0}
  $ (see Proposition \ref{propdelta}), which both have continuous trajectories since ${\bfm
  b}_i\pm{\bfm A}^n_\lambda$ have. We are then in position to
  apply Proposition \ref{prop_ito}. For any $ 0\leq    t\leq T$ and $\lambda >0$
\begin{equation*}
\begin{split}
\overline{{\bfm
  v}}_\lambda^n(Y_{t})-\overline{{\bfm
  v}}_\lambda^n(Y_{0}) & = \int_{0}^{t}[{\bfm L}\overline{{\bfm
  v}}_\lambda^n+D_t\overline{{\bfm
  v}}_\lambda^n](Y_r)\,dr+\overline{{\cal M}}^{n,\lambda}_{t} -\overline{{\cal M}}^{n,\lambda}_0  \\ & =\int_{0}^{t}[\lambda\overline{{\bfm
  v}}_\lambda^n-{\bfm
  b}_i+{\bfm A}^n_\lambda+D_t\overline{{\bfm
  v}}_\lambda^n](Y_r)\,dr+\overline{{\cal M}}^{n,\lambda}_{t} -\overline{{\cal M}}^{n,\lambda}_0  ,
  \end{split}
\end{equation*}
 where $\overline{{\cal M}}^{n,\lambda}$
   is a martingale with respect to the forward
filtration $({\cal F}_t)_{0\leq t\leq T}$, and ${\cal F}_t$ is the
$\sigma$-algebra on $\Omega $ generated by $\left\{Y_{r};0\leq r\leq
t\right\}$. From Proposition \ref{prop_ito2}, we also have
\begin{equation*}
\begin{split}
\underline{{\bfm
  v}}_\lambda^n(Y_0)-\underline{{\bfm
  v}}_\lambda^n(Y_t)& = \int_{0}^{t}[{\bfm L}^*\underline{{\bfm
  v}}_\lambda^n-D_t\underline{{\bfm
  v}}_\lambda^n](Y_r)\,dr+\underline{{\cal M}}^{n,\lambda}_{0} -\underline{{\cal M}}^{n,\lambda}_t
  \\ & = \int_{0}^{t}[\lambda\underline{{\bfm
  v}}_\lambda^n-{\bfm
  b}_i-{\bfm A}^n_\lambda-D_t\underline{{\bfm
  v}}_\lambda^n](Y_r)\,dr+\underline{{\cal M}}^{n,\lambda}_{0} -\underline{{\cal M}}^{n,\lambda}_t  ,
\end{split}
\end{equation*}
where $\underline{{\cal M}}^{n,\lambda} $ is a martingale with
respect to the backward filtration $({\cal G} _t)_{0\leq t\leq T}$,
and ${\cal G}_s$ is the $\sigma$-algebra on $\Omega $ generated by
$\left\{ Y_r;t\leq r\leq T\right\}$. Adding up these equalities, we
obtain, for any  $0\leq t\leq T$,
\begin{equation*}
\begin{split}
2\int_0^{t}{\bfm
  b}_i(Y_r)\,dr = & [\underline{{\bfm
  v}}_\lambda^n-\overline{{\bfm
  v}}_\lambda^n](Y_t)+[\overline{{\bfm
  v}}_\lambda^n-\underline{{\bfm
  v}}_\lambda^n](Y_0)+\int_{0}^{t}[\lambda(\overline{{\bfm
  v}}_\lambda^n+\underline{{\bfm
  v}}_\lambda^n)+D_t(\overline{{\bfm
  v}}_\lambda^n-\underline{{\bfm
  v}}_\lambda^n)](Y_r)\,dr\\ & +\overline{{\cal M}}^{n,\lambda}_{t} -\overline{{\cal M}}^{n,\lambda}_0+\underline{{\cal M}}^{n,\lambda}_{0} -\underline{{\cal M}}^{n,\lambda}_t
\end{split}
\end{equation*}
Fix $R>0$ and choose $\lambda=\varepsilon^2$. Integrating with
respect to the measure $\P_\pi$, we have (the sup below is taken
over $0\leq t,s \leq T$)
\begin{equation}\label{majb}
\begin{split}
\E_\pi\big[&   \sup_{|t-s|\leq \alpha}
\big|2\varepsilon\int_{s/\varepsilon^2}^{t/\varepsilon^2}{\bfm
  b}_i(Y_r)\,dr\big|\geq R\big]\\ \leq &  20R^{-2}(1+T)\varepsilon^2\big(|\underline{{\bfm
  v}}_{\varepsilon^2}^n|_2^2+|\overline{{\bfm
  v}}_{\varepsilon^2}^n|_2^2\big)+10R^{-2}T/\varepsilon^2|D_t\overline{{\bfm
  v}}_{\varepsilon^2}^n-D_t\underline{{\bfm
  v}}_{\varepsilon^2}^n|^2_2\\
+ &  5\varepsilon^2\E_\pi\big[\sup_{|t-s|\leq
\alpha}|\overline{{\cal M}}^{n,\varepsilon^2}_{t/\varepsilon^2}
-\overline{{\cal M}}^{n,\varepsilon^2}_{s/\varepsilon^2}|^2\geq
R^2\big]+5\varepsilon^2\E_\pi\big[\sup_{|t-s|\leq
\alpha}|\underline{{\cal M}}^{n,\varepsilon^2}_{s/\varepsilon^2}
-\underline{{\cal M}}^{n,\varepsilon^2}_{t/\varepsilon^2}|^2\geq
R^2\big].
\end{split}
\end{equation}
We are now going to explain how to choose $n$ to make each term of
the above right-hand side go to $0$ as $\varepsilon$ goes to $0$.\\
Since $(\lambda-{\bfm S}){\bfm w}_\lambda={\bfm
  b}_i  $ and $(\lambda-{\bfm L})\overline{{\bfm v}}^n_\lambda={\bfm b}_i-{\bfm A}^n_\lambda$,
  we can subtract these equalities and obtain, for each ${\bfm
  \varphi}\in \F$, $B^0_\lambda({\bfm w}_\lambda-\overline{{\bfm v}}^n_\lambda,{\bfm
  \varphi})={\bfm T}_H\big({\bfm w}_\lambda-{\bfm \psi}^n_\lambda,{\bfm
  \varphi}\big)$ (remind of the definition of $B^0_\lambda$ and ${\bfm T}_H$ in
  \eqref{def_dirichletform} and \eqref{def_T}). Choosing ${\bfm
  \varphi}={\bfm w}_\lambda-{\bfm \psi}^n_\lambda$, we obtain a
  first estimate
\begin{equation}\label{tight1}
 \lambda|{\bfm w}_\lambda-\overline{{\bfm v}}^n_\lambda|_2^2+(m/2)\|{\bfm w}_\lambda-\overline{{\bfm v}}^n_\lambda
 \|_1^2\leq (2m)^{-1}(C^H_1)^2\|{\bfm w}_\lambda-{\bfm
 \psi}^n_\lambda\|_1^2.
\end{equation}
Following Proposition \ref{propdelta}, we can differentiate the
equation $B^0_\lambda({\bfm w}_\lambda-\overline{{\bfm
v}}^n_\lambda,{\bfm
  \varphi})={\bfm T}_H\big({\bfm w}_\lambda-{\bfm \psi}^n_\lambda,{\bfm
  \varphi}\big)$ with respect to the time variable. So we have, for each ${\bfm
  \varphi}\in \H$, $B^0_\lambda(D_t{\bfm w}_\lambda-D_t\overline{{\bfm
v}}^n_\lambda,{\bfm
  \varphi})={\bfm T}_H\big(D_t{\bfm w}_\lambda-D_t{\bfm \psi}^n_\lambda,{\bfm
  \varphi}\big)+\partial_t{\bfm T}_H\big({\bfm w}_\lambda-{\bfm \psi}^n_\lambda,{\bfm
  \varphi}\big)-[\partial_t{\bfm T}_a+\partial_t{\bfm T}_H]\big({\bfm w}_\lambda-\overline{{\bfm v}}^n_\lambda,{\bfm
  \varphi}\big)$. Choosing ${\bfm
  \varphi}=D_t{\bfm w}_\lambda-D_t{\bfm \psi}^n_\lambda$, we obtain a
  second estimate
\begin{equation}\label{tight2}
\begin{split}
 \lambda & |D_t{\bfm w}_\lambda-D_t\overline{{\bfm v}}^n_\lambda|_2^2+(m/2)\|D_t{\bfm w}_\lambda-D_t\overline{{\bfm v}}^n_\lambda
 \|_1^2\\ & \leq (2m)^{-1}\big(C^H_1\|D_t{\bfm w}_\lambda-D_t{\bfm
 \psi}^n_\lambda\|_1+C^H_2\|{\bfm w}_\lambda-{\bfm
 \psi}^n_\lambda\|_1+(C^a_2+C^H_2)\|{\bfm w}_\lambda-\overline{{\bfm v}}^n_\lambda
 \|_1\big)^2.
 \end{split}
\end{equation}
Likewise,  \eqref{tight1} and \eqref{tight2} remain valid for  $
\underline{{\bfm v}}^n_\lambda$ instead of $\overline{{\bfm
v}}^n_\lambda$. For each fixed $\lambda>0$, we can then choose
$n_\lambda\in\nat$ large enough to ensure that $|{\bfm
w}_\lambda-\overline{{\bfm v}}^{n_\lambda}_\lambda|_2^2+\|{\bfm
w}_\lambda-\overline{{\bfm
v}}^{n_\lambda}_\lambda\|_1^2+\lambda^{-1}|D_t{\bfm
w}_\lambda-D_t\overline{{\bfm v}}^{n_\lambda}_\lambda|_2^2\leq
\lambda$ and $|{\bfm w}_\lambda-\underline{{\bfm
v}}^{n_\lambda}_\lambda|_2^2+\|{\bfm w}_\lambda-\underline{{\bfm
v}}^{n_\lambda}_\lambda\|_1^2+\lambda^{-1}|D_t{\bfm
w}_\lambda-D_t\underline{{\bfm v}}^{n_\lambda}_\lambda|_2^2\leq
\lambda$. From Proposition \ref{controlnodt}, there exists
  ${\bfm \zeta }\in (L^2(\Omega ))^d$ such that $\lambda |{\bfm w}_\lambda |_2^2+|\nabla ^\sigma {\bfm w}_\lambda- {\bfm \zeta
  }|_2\rightarrow 0$ as $\lambda$ goes to $0$. From \eqref{tight1} (with $n=n_\lambda$), $\lambda |\overline{{\bfm
v}}^{n_\lambda}_\lambda |_2^2+\lambda |\underline{{\bfm
v}}^{n_\lambda}_\lambda|_2^2\rightarrow 0$ as $\lambda$ goes to $0$.
Hence, choosing $n=n_{\varepsilon^2}$ in \eqref{majb}, all the terms
in the right-hand side except the martingale terms converge to $0$
as $\varepsilon$ goes
  to $0$.\\
  Let us now focus on the martingale terms. In order to prove the tightness
 of the two martingales, it is sufficient to prove the tightness of their brackets (see \cite{jacod} Theorem 4.13), which respectively match
$\varepsilon^2\int_0^{t/\varepsilon^2}|\nabla^\sigma\overline{{\bfm
  v}}_{\varepsilon^2}^{n_{\varepsilon^2}}(Y_r)|^2\,dr $ and $\varepsilon^2\int_0^{t/\varepsilon^2}|\nabla^\sigma\underline{{\bfm
  v}}_{\varepsilon^2}^{n_{\varepsilon^2}}(Y_r)|^2\,dr $. Note that
  $|\nabla^\sigma\overline{{\bfm
  v}}_{\varepsilon^2}^{n_{\varepsilon^2}}- {\bfm \zeta
  }|_2\rightarrow 0$ as $\varepsilon$ tends to $0$ so that the process $t\mapsto\varepsilon^2\int_0^{t/\varepsilon^2}|\nabla^\sigma\overline{{\bfm
  v}}_{\varepsilon^2}^{n_{\varepsilon^2}}(Y_r)|^2\,dr $ has the
  same limit in $C([0,T];\R)$ as the process $t\mapsto\varepsilon^2\int_0^{t/\varepsilon^2}|{\bfm
  \zeta}(Y_r)|^2\,dr $. Finally, for each fixed $t$, Theorem \ref{theo_erg} proves that  $\varepsilon^2\int_0^{t/\varepsilon^2}|{\bfm
  \zeta}(Y_r)|^2\,dr $ converges to the deterministic non-decreasing process $t\int_\Omega
|{\bfm \zeta }|_2^2 \,d\pi $ in $L^1$ under the measure $\P_\pi$.
Then Theorem 3.37 in \cite{jacod} says that the brackets are tight
in $C([0,T];\reel)$. The same arguments remain valid for the
brackets of $\underline{{\cal
M}}^{n_{\varepsilon^2},\varepsilon^2}$. Hence, the right-hand side
in \eqref{majb} converges to $0$ as $\varepsilon$ goes to $0$ and
the tightness of $t\mapsto \varepsilon X^\omega_{t/\varepsilon^2}$
follows. \qed

%%%%%%%%%%%%%%%%%%%%%%%%%%%%%%%%%%%%%%%%%%%%%%%%%%%%%%%%%%%%%%%%%%%%%%%%%%%%%%%%%%%%%%%%%%


\begin{thebibliography}{60cm}

%%%%%%%%%%%%%%%%%%%%%%%%%%%%%%%%%%%%%%%%%%%%%%%%%%%%%%%%%%%%%%%%%%%%%%%%%%%%%%%%%%%%%%%%%%

\bibitem{bhat}
\textsc{Bhattacharya, Gupta, Walker}, Asymptotics of solute
dispersion in periodic media, SIAM J. Appl. Math. 49 (1989), n° 1,
86-98.
%\bibitem{pianitski}
%A. Bourgeat \& A. Pianitski, Approximations of effective coefficients in stochastic
%homogenization, Annales de l'institut Henri Poincarré, Probabilités et statistiques, PR (2004), 153-165.
\bibitem{fannjang}
\textsc{Fannjiang}, \textsc{Komorowski}, An invariance principle for
diffusion in turbulence, Annals of Probability, 1999, vol. 27, No.
2, 751-781.
\bibitem{fannjang4}
\textsc{Fannjiang}, \textsc{Komorowski}, Diffusion approximation for
particle convection in Markovian Flows, Bull. Polish Acad. Sci.
Math., 2000, vol. 48, No. 3, 253-275.
\bibitem{fannjang5}
 \textsc{Fannjiang}, \textsc{Komorowski}, Invariance principle for a diffusion in a Markov field,
Bull. Polish Acad. Sci. Math., 2001, vol. 49, No. 1, 45-65.
%\bibitem{delarue}
%F. Delarue, thèse: {\it EDSPR. Application à l'homogénéisation des EDP quasi-linéaires},
%Université de Provence (Marseille, France), janvier 2002.
\bibitem{fukushima}
\textsc{Fukushima} , \textsc{Oshima, Takeda}, Dirichlet Forms and
Symmetric Markov Processes, De Gruyter Studies in Mathematics 19,
Walter de Gruyter, Berlin and Hawthorne, New York, 1994.
%\bibitem{ggoss}
%\textsc{Gantert, Garnier, Olla, Shi, Sznitman}, Milieux aléatoires,
%édité par Francis Comets et Etienne Pardoux, Panoramas et Synthèses,
%numéro 12, Société mathématique de France, 2001, pp. 75-99.
\bibitem{jacod}
\textsc{Jacod, Shiryaev}, Limit Theorems for Stochastic Processes,
Grundlehren der mathematischen Wissenschaft 288, Springer-Verlag.
\bibitem{kipnis}
\textsc{Kipnis , Varadhan }, Central limit theorem for additive
functionals of reversible Markov processes and application to simple
exclusion, Ann. Probab.  28 (2000), no. 1 , 277-302.
\bibitem{kozlov}
\textsc{Kozlov}, The Method of Averaging and Walks in Inhomogeneous
Environments, Russian Math. Surveys. (1985), 40, 73-145.
\bibitem{kusuoka}
\textsc{Kusuoka, Stroock}, Long time estimates for the heat kernel
associated with a uniformly subelliptic symmetric second order
operator, Annals of Mathematics, 127 (1988), 165-189.
%\bibitem{lejay}
%A.~Lejay, Méthodes probabilistes pour l'homogénéisation des
%opérateurs sous forme divergence: cas linéaire et semilinéaires,
%thèse, Marseille, 2000.
\bibitem{olla1}
\textsc{Landim, Olla, Yau}, Convection-diffusion equation with
space-time ergodic random
 flow, Probability theory and related fields 112 (1998), 203-220.
\bibitem{olla2}
\textsc{Komorowski, Olla}, On homogenization of time-dependent
random flows, Probability
 theory and related fields 121 (2001), 98-116.
 \bibitem{ma}
\textsc{Ma, Röckner}, Introduction to the Theory of (Non-Symmetric)
Dirichlet Forms, Universitext, Berlin Heidelberg, Springer-Verlag,
1992.
%\bibitem{michel}
%D.Michel, E. Pardoux, {\it An introduction to Malliavin Calculus and
%some of its application},.
\bibitem{oelschlager}
\textsc{Oelschläger},  Homogenization of a diffusion process in a
divergence free random field, Annals of Probability, 1988, 16,
1084-1126.
\bibitem{olla}
\textsc{Olla}, Homogenization of diffusion processes in Random
Fields, Cours de l'école doctorale, Ecole polytechnique, 1994.
\bibitem{osada}
\textsc{Osada}, Homogenization of diffusion with random stationary
coefficients, Lecture Notes in Math. 1021 (1982), pp. 507-517.
\bibitem{pardouxdeg}
\textsc{Pardoux}, Homogenization of periodic linear degenerate PDEs,
LATP, Université de Provence, Marseille, 2005.
%\bibitem{pardouxtopics}
%E. Pardoux, {\it Topics in the probabilistic approach to
%homogenization.}
\bibitem{pardouxper}
\textsc{Pardoux}, Homogenization of linear and semilinear second
order parabolic PDEs with periodic coefficients: a probabilistic
approach, J. Funct. Anal., 167, 498-520.
\bibitem{protter}
\textsc{Protter} P.,  Stochastic integration and differential
equations, A new approach, Applications of mathematics,
Springer-Verlag, Berlin, 1990.
\bibitem{revuz}
\textsc{Rudin},  Analyse fonctionnelle, Ediscience International,
1995.
\bibitem{sethu}
\textsc{Sethuraman, Varadhan, Yau}, Diffusive limit of a tagged
particle in asymmetric simple exclusion processes, Commun. Pure and
Appl. Math. (2000), 53, 972-1006.
\bibitem{sznitman}
\textsc{Sznitman, Zeitouni}, An invariance principle for isotropic
diffusions in random environment, C.R. Acad. Sci. Paris, Ser. I 339
(2004), 429-434.
\bibitem{wu}
\textsc{Wu}, Forward-Backward martingale decomposition and
compactness results for additive functionals of stationary ergodic
Markov processes, Ann. Inst. Henri Poincaré 35 (1999), 121-141.

\end{thebibliography}
\end{document}